\documentclass{nm}
\renewcommand\thisnumber{x}
\renewcommand\thisyear {202x}

\renewcommand\thisvolume{xx}
\usepackage{charter}
\usepackage[charter]{mathdesign}

\makeatletter
\let\algorithm\@undefined
\let\endalgorithm\@undefined
\makeatother

\usepackage[linesnumbered,ruled,vlined]{algorithm2e}
 \usepackage{amsfonts}
\usepackage{amsmath}
\usepackage{array}
\usepackage{cases}
\usepackage{color}
\usepackage{graphicx}
\usepackage{subfigure}

\newcommand{\lvertiii}{{\left\vert\kern-0.25ex\left\vert\kern-0.25ex\left\vert}}
\newcommand{\rvertiii}{{\right\vert\kern-0.25ex\right\vert\kern-0.25ex\right\vert}}
\newtheorem{assume}[theorem]{Assumption} %asuume用 theorem的计数器
\newcommand{\rv}[1]{#1}%{{\color{red} #1}}

\begin{document}

%%%%% title and author(s):
% \markboth{Author(s)}{Short Title}
% \title{Title}

\markboth{X. Dai, L. Zhang and A. Zhou}{Orthogonality preserving schemes for approximating the KS orbitals}
\title{Convergent and orthogonality preserving schemes for approximating the Kohn-Sham orbitals
	\thanks{This work was supported by the National Key R $\&$ D Program of China under grants
		2019YFA0709600 and 2019YFA0709601, the National Natural Science Foundation of China under
		grants 12021001.}}

% single author:
% \author[AUTHOR]{AUTHOR\corrauth}
% \address{address of AUTHOR}
% \email{{\tt email address of AUTHOR} (AUTHOR)}

% multiple authors:
% Please mark \corrauth after the name of the corresponding author.
% different addresses:
%\author[AUTHOR1 and AUTHOR2]{AUTHOR1\affil{1}\comma\corrauth and AUTHOR2\affil{2}}
%\address{\affilnum{1}\ address of AUTHOR1\\
%\affilnum{2}\ address of AUTHOR2}
%
%same address:
%\author[AUTHOR1, AUTHOR2 and AUTHOR3]{AUTHOR1, AUTHOR2\corrauth and AUTHOR3}
%\address{address of AUTHOR1, AUTHOR2 and AUTHOR3}
%
%\emails{{\tt email of AUTHOR1} (AUTHOR1), {\tt email of AUTHOR2} (AUTHOR2), {\tt email of AUTHOR3} (AUTHOR3)}
%

\author[Xiaoying Dai, Liwei Zhang and Aihui Zhou]{Xiaoying Dai\corrauth, Liwei Zhang and Aihui Zhou}
\address{LSEC, Institute of Computational Mathematics and Scientific/Engineering Computing,
	Academy of Mathematics and Systems Science, Chinese Academy of Sciences,  Beijing 100190, China; and School of Mathematical Sciences,
	University of Chinese Academy of Sciences, Beijing 100049, China.}

\emails{{\tt daixy@lsec.cc.ac.cn} (Xiaoying Dai), {\tt zhanglw@lsec.cc.ac.cn} (Liwei Zhang), {\tt azhou@lsec.cc.ac.cn} (Aihui Zhou)}

\begin{abstract}
	To obtain convergent numerical approximations without using any orthogonalization operations is of great importance in electronic structure calculations. In this paper, we propose and analyze a class of  iteration schemes for the discretized Kohn-Sham Density Functional Theory model, with which the iterative approximations are guaranteed to converge to the Kohn-Sham orbitals exponentially without any orthogonalization as long as the initial orbitals are orthogonal and the time step sizes are given properly. 
	In addition,  we present a feasible and efficient approach to get suitable time step sizes \rv{and report some numerical experiments to validate our theory.}
\end{abstract}


%%%% AMS subject classifications %%%%

%%%% maketitle %%%%% 
\maketitle

\section{Introduction}\label{sec-intro}
Electronic structure calculations play an important role in numerous fields such as quantum chemistry, materials science and drug design. Due to the good balance of accuracy and computational cost, the Kohn-Sham Density Function Theory (DFT) model \cite{HK,KS65,Matin04,parr-yang94,payne92} \rv{has became} one of the most widely used models in electronic structure  calculations which   is usually treated as either a nonlinear eigenvalue problem (Kohn-Sham equation) or an orthogonality constraint minimization problem (Kohn-Sham total energy  direct minimization problem). %DWZ

In the literature,  there are a number of works on the design and analysis of numerical methods for solving the Kohn-Sham equation (see, e.g., \cite{cances12,chen13,chen-dai14,dai-zhou15,lin-lu-ying2019,Saad10,zhang-shen08}  and references cited therein). 
To obtain the solution of this nonlinear eigenvalue problem, we observe that some self consistent field (SCF) iterations are usually used \cite{Matin04} (see also \cite{anderson65,Cances01,johnson88,lin-yang13,pulay80,pulay82,zhou-wang-2018}). Unfortunately, the convergence of SCF iterations is uncertain. We understand that its convergence has indeed been investigated when there is a \rv{sufficiently} large gap between the occupied and unoccupied states and the second-order derivatives of the exchange correlation functional are uniformly bounded from above \cite{BLL,LWWUY,LWWY,YGM}, which is  important in the theoretical point of view.
It becomes significant to investigate the convergence of SCF iterations when the gap is not
large in application.

We see that an alternative approach   to obtain the ground states is to solve  the Kohn-Sham total energy minimization  problem, which is an orthogonality constrained minimization problem \cite{payne92}. The direct minimization approach attracts the attention of many researchers in recent years \cite{Cances21,Frey09,Marzari97}, and many different kinds of optimization methods  are applied to electronic structure calculations and investigated (see, e.g., \cite{DLZZ,DLZxZ,GLCY,GLY,schneider09,yang-wang07,ZBJ,ZZWZ}). 
   
For solving either the nonlinear eigenvalue problem or the orhtogonality constrained minimization problem,   except for  few works such as \cite{GLY},  the orthogonalization procedure is   usually required, which is very expensive and limits the parallel scalability in numerical implementation.
 
Recently, Dai et al. proposed a gradient flow based Kohn-Sham DFT model \cite{DWZ} that is a time evolution problem and is completely different from either the nonlinear eigenvalue problem or the orthogonality constrained minimization problem. It is proved in \cite{DWZ} that the flow   of the new model is orthogonality preserving, and \rv{the solution can evolve} to the ground state. Consequently, the gradient
flow based model provides a novel and attractive approach for solving Kohn-Sham
DFT apart from the eigenvalue problem model and the energy minimization
model. \rv{In other words,} the gradient flow based model is quite promising
in \rv{ground state} electronic structure calculations and deserves further investigation.
\rv{For the sake of clarity, we would like to mention that the gradient flow based Kohn-Sham DFT model is different from the time dependent Kohn-Sham equation in \cite{M3GR, RG1984, YSHH}.}

In this paper,  we  propose a general framework of orthogonality preserving schemes  that produce efficient approximations of the Kohn-Sham orbitals with the help of the gradient flow based model.  In addition, we prove the global convergence and local convergence
rate of the new schemes under some mild assumptions. We also provide some typical 
choices for the auxiliary mapping appeared in the framework, and a feasible and efficient
approach to obtain the desired time step sizes that satisfy the assumptions required in
the  analysis, which result in several typical orthogonality preserving schemes
that can produce convergent approximations of the Kohn-Sham orbitals.

The rest of the paper is organized as follows. In Section \ref{sec-mod}, we briefly review the gradient flow based Kohn-Sham DFT model and some notation frequently used \rv{throughout} this paper.  We then propose a framework \rv{for} orthogonality preserving schemes for solving the discretized Kohn-Sham model in Section \ref{sec-scheme} and prove its global convergence as well as local convergence rate under some reasonable assumptions and with proper time step size.  Then, in  Section \ref{sec-ii}, we  provide some specific choices for the  auxiliary mapping  and the time step size. \rv{We then report some numerical results obtained by the proposed schemes in Section \ref{sec-ne} to verify our theory.} Finally, we give some concluding remarks in Section \ref{sec-cln}. 

\section{Kohn-Sham DFT Models}\label{sec-mod}
\subsection{  Classical Kohn-Sham DFT model}
According to Kohn-Sham density functional theory \cite{KS65}, the ground state of a system can be obtained by solving
\begin{equation}\label{dis-emin}
\begin{split}
&\inf_{U\in {\big(H^1(\mathbb{R}^3)\big)}^N} \ \ \ E_{\textup{KS}}(U) \\
& s.t.\ \   U^TU = I_N,
\end{split}
\end{equation}
where $U = (u_1, \dots, u_N)\in {\big(H^1(\mathbb{R}^3)\big)}^N$, $U^TV = \big(\langle u_i,v_j\rangle_{L^2(\mathbb{R}^3)}\big)_{i,j=1}^N, \forall U,V\in\big(H^1(\mathbb{R}^3)\big)^N$,  and the objective functional $E_{\textup{KS}}(U)$ reads as
\begin{eqnarray}\label{energy}
E_{\textup{KS}}(U)&=&\frac{1}{2} \int_{\mathbb{R}^3} \sum_{i=1}^N|\nabla u_i(r)|^2 dr +
\frac{1}{2}\int_{\mathbb{R}^3}\int_{\mathbb{R}^3}\frac{\rho(r)\rho(r')}{|r-r'|}drdr' \nonumber\\
&&+\int_{\mathbb{R}^3} V_{ext}(r)\rho(r)dr+ \int_{\mathbb{R}^3}
\varepsilon_{xc}(\rho)(r)\rho(r)dr.
\end{eqnarray}
Here, $N$ denotes the number of electrons, $\{u_i\}_{i=1,2,\cdots,N}$ are usually called the Kohn-Sham orbitals,
$\rho(r)=\sum \limits_{i=1}^N|u_i(r)|^2$ is the electronic density (we assume each Kohn-Sham orbital is occupied by one electron here),
$V_{ext}(r)$ is the external potential generated by the nuclei, and $\varepsilon_{xc}(\rho)(r)$ is the exchange-correlation functional
which is not known explicitly. In practice, some approximation such as local density approximation (LDA),
generalized gradient approximation (GGA) or some other approximations has to be used \cite{Matin04}.

We see that the feasible set of \eqref{dis-emin} is a Stiefel manifold which is defined as
\begin{equation}
\mathcal{M}^N = \{U\in {\big(H^1(\mathbb{R}^3)\big)}^N: U^TU = I_N\}.
\end{equation}
To get rid of the nonuniqueness of the minimizer caused by the invariance of the energy functional under orthogonal transformations to the Kohn-Sham orbitals (i.e., $E(U) = E(UP), \iffalse~\rho(U)=\rho(UP), \fi~\forall P\in\mathcal{O}^N$ with $\mathcal{O}^N$ being the set of orthogonal matrices of order $N$), we, following \cite{DLZZ,DWZ,DZZ}, consider \eqref{dis-emin} on the Grassmann manifold $\mathcal{G}^N$ which is a quotient manifold of Stifel manifold, that is
\begin{equation*}
\mathcal{G}^N = \mathcal{M}^N/\sim.
\end{equation*}
Here, $\sim$ denotes the equivalence relation which is defined as: $\hat{U} \sim U$,  if and only if there exists $P \in \mathcal{O}^N$, such that $\hat{U} = UP$. For any $U \in \mathcal{M}^N$, we denote
\begin{equation*}
[U] = \{U P: P \in  \mathcal{O}^N\},
\end{equation*}
then the Grassmann manifold $\mathcal{G}^N$ can be formulated as
\begin{equation*}
\mathcal{G}^N= \{[U]: U\in\mathcal{M}^N\}.
\end{equation*}
For $[U] \in \mathcal{G}^N$, the tangent space of $[U]$ on $\mathcal{G}^N$ is the following set 
\begin{equation}\label{tan}
\mathcal{T}_{[U]}\mathcal{G}^N = \{W\in V^N:
W^TU   = {\bf 0} \in \mathbb{R}^{N\times N}\}.
\end{equation}

In this paper, we assume that \eqref{dis-emin} achieves its minimum in $\mathcal{G}^N$, which implies that
\eqref{dis-emin} is equivalent to
\begin{equation}\label{Emin-Grass}
\min_{[U] \in \mathcal{G}^N} \ \ \ E(U).
\end{equation}

In addition, we see from \cite{AbMaSe} that the Grassmann gradient of $E_{\textup{KS}}(U)$ is
\begin{equation*}
\nabla_G E_{\textup{KS}}(U)  = \nabla E_{\textup{KS}}(U) - UU^T\nabla E_{\textup{KS}}(U),
\end{equation*}
where $$\nabla E_{\textup{KS}}(U) = \mathcal{H}(\rho)U, ~ \forall U\in{\big(H^1(\mathbb{R}^3)\big)}^N$$ is the Euclidean gradient of $E_{\textup{KS}}(U)$,
$$\mathcal{H}(\rho) = -\frac{1}{2}\Delta + V_{ext} + \int_{\mathbb{R}^3}\frac{\rho(r')}{|r-r'|}dr' +
v_{xc}(\rho),$$
is symmetric and $$v_{xc}(\rho)
= \frac{\delta(\rho\varepsilon_{xc}(\rho))}{\delta\rho}.$$

For any $U\in\mathcal{M}^N$, we see from \cite{DWZ} that
 $\nabla_G E_{\textup{KS}}(U)  = \mathcal{A}_UU,$
 where $$\mathcal{A}_U = \nabla E_{\textup{KS}}(U)U^T - U\nabla E_{\textup{KS}}(U)^T$$ is anti-symmetric, \rv{i.e., ${\mathcal{A}_U}^T = -\mathcal{A}_U$}.
Furthermore, the Hessian of $E_{\textup{KS}}(U)$ on $\mathcal{G}_N$ has the form \cite{DLZZ}
\begin{equation*}
\begin{split}
\nabla^2_GE(U)[D_1, D_2] &= \text{tr} ( {D_1}^T\mathcal{H}(\rho)
D_2) - \text{tr}( {D_1}^TD_2U^T\mathcal{H}(\rho)U) \\
&+ 2\int_{\mathbb{R}^3}\int_{\mathbb{R}^3}\frac{(\sum_iu_i(r)d_{1,i}(r))(\sum_ju_j(r')d_{2,j}(r'))}
{|r-r'|}drdr' \\
&+ 2\int_{\mathbb{R}^3}\frac{\delta^2(\varepsilon_{xc}(\rho)\rho)}{\delta\rho^2}(r)
(\sum_iu_i(r)d_{1,i}(r))(\sum_ju_j(r)d_{2,j}(r))dr
\end{split}
\end{equation*}
provided that the total energy functional is \rv{twice} differentiable, or more specifically,
the approximated exchange-correlation functional is \rv{twice} differentiable.
Here, 
\[
D_i = (d_{i,1}, d_{i,2}, \cdots, d_{i,N}) \in \mathcal{T}_{[U]}\mathcal{G}^N,~ i = 1,2.
\]

The ground state of a system can also be obtained by considering the Euler-Lagrange equation of \eqref{dis-emin}, which reads as
\begin{equation}\label{dis-evp}
\begin{cases}
\mathcal{H}(\rho)U = U\Lambda, \\
U\in\mathcal{M}^N.
\end{cases}
\end{equation}
The nonlinear eigenvalue problem \eqref{dis-evp} is indeed the so-called Kohn-Sham equation. For decades, the Kohn-Sham DFT models are investigated as either a minimization problem \eqref{Emin-Grass} or an eigenvalue problem \eqref{dis-evp}. 

In practice, we may \rv{discretize} the Kohn-Sham energy minimization model \eqref{Emin-Grass} as well as the nonlinear eigenvalue model \eqref{dis-evp} by, e.g., the plane wave method, the local basis set method, or real space methods. More details about the discretization methods can be found in, for instance, the review paper \cite{Saad10}. If we choose a $N_g$-dimension space $V_{N_g}\subset H^1(\mathbb{R}^3)$
to approximate $H^1(\mathbb{R}^3)$, then the associated discretized Kohn-Sham model
can be formulated as
\begin{equation}\label{E_dKS}
\min_{[U]\in \mathcal{G}^N_{N_g}} \ \ \ E_{KS}(U),
\end{equation}
or
\begin{equation}\label{EVP_dKS}
\begin{cases}
\mathcal{H}(\rho)U = U\Lambda, \\
U\in\mathcal{M}_{N_g}^N,
\end{cases}
\end{equation}
where $\mathcal{G}^N_{N_g}$ is the discretized Grassmann manifold defined by
\begin{equation*}
\mathcal{G}^N_{N_g} = \mathcal{M}^N_{N_g}/\sim,
\end{equation*}
and
\begin{equation*}
\mathcal{M}^N_{N_g} = \{U \in(V_{N_g})^N:   U^TU  = I_N\}
\end{equation*}
is the discretized Stiefel manifold with the \rv{equivalence} relation $\sim$ having the similar meaning to what we have mentioned. Usually, $N_g\gg N$. We should point out that the operators on the discretized manifold, such as the Grassmann gradient and the Grassmann Hessian, have exactly the same forms as those on the continuous manifold. 

\subsection{Gradient flow based  Kohn-Sham DFT model}
Different from the minimization model \eqref{E_dKS} and the eigenvalue model \eqref{EVP_dKS}, \rv{a gradient flow based Kohn-Sham DFT model was proposed in \cite{DWZ},
%\rv{in which the orbitals $U$ is considered to be evolved with respect to the time variable $t$, whose evolution is governed by the following equation:}
which has the following form:}
\begin{equation}\label{Gra-Flow}
\begin{cases}
\frac{d}{dt}U(t) = -\nabla_G E(U(t)), \forall t\in\mathbb{R}^+, \\
U(0) = U_0\in\mathcal{M}^N_{N_g}.
\end{cases}
\end{equation}
\rv{Here, as $U_0$ is required to be orthogonal, we see from} \cite{DWZ} that for the gradient flow based model \eqref{Gra-Flow}, there hold
\begin{equation}\label{ortho-pre}
U(t)\in\mathcal{M}^N_{N_g}, \forall t\ge0,
\end{equation}
and
\begin{equation}\label{non-increase}
\frac{d}{dt}E(U(t)) = -\|\nabla_GE(U(t))\|^2\le0, \forall t\ge0.
\end{equation}
Besides, it is  proved in  \cite{DWZ} that  the norm of the extended gradient of energy functional   exponentially decays to zero over time $t$, and the \rv{solution $U(t)$} will evolve to the ground state under some mild assumptions. \rv{We mention that the detailed derivation of the gradient flow based Kohn-Sham model can be found in \cite{DWZ}.}
 
 \section{A general framework of orthogonality preserving schemes}\label{sec-scheme}
 With the help of the gradient flow based model (2.9),   we design a general \rv{framework} that enables  us to obtain a  class of orthogonality preserving schemes for getting  convergent approximations of Kohn-Sham orbitals. To propose our numerical schemes, we first introduce the partition of the time interval
$$0 = t_0<t_1<\cdots<t_n<\cdots.$$
 In addition, we denote $\Delta t_n = t_{n+1} - t_n$ and use $U_n$ to symbolize the approximation of $U(t_n)$ for $n\in\mathbb{N}_0$ where $\mathbb{N}_0$ is the set of nonnegative integers.

\subsection{Scheme framework}
Given $U_0\in\mathcal{M}_{N_g}^N$, we consider the following recursive scheme on the interval $[t_n,t_{n+1})$:
\begin{equation}\label{generate}
\begin{cases}
 \tilde{U}(t)-U_n = -(t-t_n)\mathcal{A}_{U^{\textup{Aux}}(t)}\frac{U_n+\tilde{U}(t)}{2}, ~ t\in[t_n,t_{n+1}), \\
U_{n+1}=\tilde{U}(t_{n+1}^-).
\end{cases}
\end{equation}
Here $U^{\textup{Aux}}:\mathbb{R}\to(V_{N_g})^N$ is a piecewise smooth auxiliary mapping which satisfies \rv{the interpolation condition} $U^{\textup{Aux}}(t_n) = U_n$  for all $n$. We hence name our schemes as interpolation based schemes. \rv{Note that this interpolation condition is the only restriction on $U^{\text{Aux}}$, which makes our framework quite flexible.} % The numerical impact of $U^{\text{Aux}}$ requires further studies and is somewhat beyond the scope of this paper.
%In general, we are not able to determine $U^{\textup{Aux}}(t)$ over the whole time interval in advance since $\{U_{k}\}_{k\ge n+1}$ are unknown at the $n$-th time interval $[t_n,t_{n+1})$. We may instead define $U^{\textup{Aux}}(t)$ over each time subinterval successively.
As a result, we obtain the following \rv{framework for} interpolation based scheme (Algorithm \ref{A1}) for solving \eqref{Gra-Flow}.

\begin{algorithm}\label{A1}
	\caption{A framework for interpolation based scheme}
	Given $\epsilon>0$,
	initial orbitals $U_0\in\mathcal{M}_{N_g}^N$,
	calculate the gradient $\nabla_G E(U_0)$ and let $n = 0$, \rv{$t_0=0$}\;
	\While{$\| \nabla_G E(U_n) \|> \epsilon$}{
		Choose a suitable $\Delta t_n>0$ and let $t_{n+1}=t_n+\Delta t_n$;
		
		Define $U^{\textup{Aux}}(t), t\in [t_n,t_{n+1})$ such that $U^{\textup{Aux}}(t_n)=U_n$;
		
		 		Update $U_{n+1} = \lim_{t\to t_{n+1}^-}\tilde{U}(t)$ with $\tilde{U}(t)$ satisfying 
		\begin{equation}\label{inner}
		\tilde{U}(t)-U_n = -(t-t_n)\mathcal{A}_{U^{\textup{Aux}}(t)}\frac{U_n+\tilde{U}(t)}{2}, ~ t\in [t_n,t_{n+1});
		\end{equation}
		
		Let $n=n+1$, calculate the gradient $\nabla_G E(U_n)$\;
	}
\end{algorithm}
  
We see from Algorithm \ref{A1} that the time step sizes in our scheme can be provided step by step and adaptively, i.e., we may make a full use of the information obtained during the iteration to determine a suitable time step size at each step. 
Besides, we point out here that if the definition of $U^{\textup{Aux}}(t),~ t\in [t_n,t_{n+1})$ is independent of $\tilde{U}(t),~ t\in[t_n, t_{n+1})$ at the $n$-th iteration, then the problem \eqref{inner} is said to be linear, and the corresponding scheme is called an explicit scheme. Otherwise, it is an implicit scheme. The following Theorem \ref{thm-orthopreserve} shows that \rv{Algorithm \ref{A1}} preserves the orthogonality of iterations no matter whether it is explicit or not.

\begin{theorem}\label{thm-orthopreserve}
	 If $\{U_n\}_{n\in\mathbb{N}_0}$ is produced by Algorithm \ref{A1}, then $\{U_n\}_{n\in\mathbb{N}_0}\subset \mathcal{M}_{N_g}^N$. 	
\end{theorem}

\begin{proof}
	By rearranging \eqref{generate}, we have that
	\begin{equation}\label{update-formula}
	U_{n+1} = \big(I+\frac{\Delta t_n}{2}\mathcal{A}_{U^{\textup{Aux}}(t_{n+1}^-)}\big)^{-1}\big(I-\frac{\Delta t_n}{2}\mathcal{A}_{U^{\textup{Aux}}(t_{n+1}^-)}\big) U_n.
	\end{equation}
	Since $\mathcal{A}_{U^{\textup{Aux}}(t_{n+1}^-)}$ is anti-symmetric, we see that $$\big(I+\frac{\Delta t_n}{2}\mathcal{A}_{U^{\textup{Aux}}(t_{n+1}^-)}\big)^{-1}\big(I-\frac{\Delta t_n}{2}\mathcal{A}_{U^{\textup{Aux}}(t_{n+1}^-)}\big)$$
	forms a Cayley transformation. Hence, $U_{n+1}\in\mathcal{M}_{N_g}^N$ as long as $U_n\in\mathcal{M}_{N_g}^N$. Note that $U_0\in\mathcal{M}_{N_g}^N$, we complete the \rv{proof by induction.}
\end{proof}

In fact, we see from the proof of Theorem \ref{thm-orthopreserve} that $\tilde{U}(t)$ is orthogonal for all $t\rv{\in\mathbb{R}^+}$, namely, $\tilde{U}(t)\rv{\subset}\mathcal{M}_{N_g}^N$. In addition, we see that for any explicit scheme, the orbitals $U_{n+1}$ can be updated simply by \eqref{update-formula}.
 Therefore, to update $U_{n+1}$ at each iteration of an explicit scheme, the main cost is to compute the inverse of $$I+\frac{\Delta t_n}{2}\mathcal{A}_{U^{\textup{Aux}}(t_{n+1}^-)},$$ which is a $N_g$-dimensional matrix inverse problem and is very expensive \rv{to obtain}.  Even though we can deal with it by solving the corresponding linear system using some iterative methods, it is still not cheap, especially when $N_g$ is large.

Fortunately, we observe that $\mathcal{A}_U~(\forall U\in V^N)$ is \rv{anti-symmetric} and has the following factorization
\begin{equation*}       
\mathcal{A}_U =
\left(                 
\begin{array}{cc}   
\nabla E(U) & U  
\end{array}
\right)   
\left(      
\begin{array}{c}   
U^T \\  
-\nabla E(U)^T
\end{array}
\right).
\end{equation*}
 Hence, by applying the Sherman-Morrison-Woodbury (SMW) formulae \cite{DWZ, DingZhou}, we have
\begin{eqnarray}\label{SMW}
&&\big(I+\frac{\Delta t_n}{2}\mathcal{A}_{U}\big)^{-1} = I - \frac{\Delta t_n}{2}
\left(                 
\begin{array}{cc}   
\nabla E(U) & U  
\end{array}
\right)   \notag \\
&&\big[
I_{2N}+\frac{\Delta t_n}{2} 
\left(                 
\begin{array}{cc}   
U^T\nabla E(U) & U^TU\\
-\nabla E(U)^T\nabla E(U) & -\nabla E(U)^TU
\end{array}
\right)   
\big]^{-1}
\left(      
\begin{array}{c}   
U^T \\  
-\nabla E(U)^T
\end{array}
\right),
\end{eqnarray}
which reduces the dimension of the matrix inverse problem significantly from $N_g$ to $2N$. Therefore, we only need to deal with linear system of dimension $2N$.

\subsection{Numerical analysis} 
Note that \eqref{non-increase} indicates the energy functional is non-increasing with respect to $t$, we may impose the following assumption on the time step sizes to maintain \rv{a} similar property. We will show the existence of the desired time partition and introduce an efficient strategy to obtain such time step sizes in the next section.

\begin{assume}\label{Assume-energy-decrease}
	The sequence $\{t_n\}_{n=0}^{\infty}$ satisfies
	\begin{equation}\label{sum_unbound}
	\sum_{n=0}^{\infty}\Delta t_n = +\infty,~ \textup{i.e.,}~ \displaystyle\lim_{n\to\infty}t_n = +\infty
	\end{equation}
	and
	\begin{align}\label{back-cond}
	E(U_{n+1})-E(U_n) = E(\tilde{U}(t_{n+1}^-)) - E(U_n)\leq -\eta \Delta t_n\|\nabla_GE(U_n)\|^2,  n\in\mathbb{N}_0
	\end{align}
with $\eta>0$ being a given parameter.
\end{assume}
 
The condition \eqref{sum_unbound} on Assumption \ref{Assume-energy-decrease} is simple   and reasonable, as we are discretizing an infinite time period. Meanwhile, the condition \eqref{back-cond} follows from \eqref{non-increase}, which 
indicates that the finite difference approximation of the temporal derivative stated in the left hand side of \eqref{non-increase} is somewhat comparable to $\|\nabla_G E(U(t))\|^2$.

Under Assumption \ref{Assume-energy-decrease}, we obtain the following asymptotic behaviour for \rv{the} approximate\rv{d} solution of \eqref{Gra-Flow} produced by Algorithm \ref{A1}.
\begin{theorem}\label{convergence}
	If
	 the sequence $\{t_n\}_{n\in\mathbb{N}_0}$ satisfies Assumption \ref{Assume-energy-decrease}, then for the sequence $\{U_n\}_{n\in\mathbb{N}_0}$ produced by Algorithm \ref{A1} with initial guess $U_0\in\mathcal{M}_{N_g}^N$, there holds
	$$\displaystyle \liminf_{n\to\infty}\|\nabla_G E(U_n)\|=0.$$
\end{theorem}
\begin{proof}
	 We see from Assumption \ref{Assume-energy-decrease} that
	$$E(U_n) - E(U_{n+1})\ge\eta\Delta t_n\|\nabla_G E(U_n)\|^2.$$
	Hence,
	\begin{eqnarray*}
		\displaystyle
		E(U_0)-E_{\min}\ge\sum_{n=0}^\infty \big(E(U_n) - E(U_{n+1})\big) \\
		\ge\eta\sum_{n=0}^\infty\Delta t_n\|\nabla_G E(U_n)\|^2,
	\end{eqnarray*}
	where $E_{\min}$ is the minimum of the energy functional $E(U)$.
	Thus, 
	\begin{equation}\label{sumunbound}
	\sum_{n=0}^\infty\Delta t_n\|\nabla_G E(U_n)\|^2<\infty.
	\end{equation}
	If $\displaystyle \liminf_{n\to\infty} \|\nabla_G E(U_n)\|>0$, then there exists $\epsilon_0>0$ such that
	$$\|\nabla_G E(U_n)\|\ge\epsilon_0, \forall n\in\mathbb{N}_0.$$
	Hence, we obtain from \eqref{sum_unbound} that
	$$\sum_{n=0}^\infty\Delta t_n\|\nabla_G E(U_n)\|^2\ge\epsilon_0^2\sum_{n=0}^\infty\Delta t_n=\infty,$$
	which is \rv{contradictory} to \eqref{sumunbound}. As a result,
	$$\liminf_{n\to\infty} \|\nabla_G E(U_n)\|=0.$$
\end{proof}

We see that a sufficient condition for \eqref{sum_unbound} is $\Delta t_n>\tau, \forall n\in\mathbb{N}_0$ for some $\tau>0$. Under this setting, Theorem \ref{convergence} indicates that the sequence $\{U_n\}_{n=0}^\infty$ produced by Algorithm \ref{A1} will converge to \rv{an} equilibrium point of \eqref{Gra-Flow} (at least for a subsequence). If the equilibrium point \rv{(denoted by $U^\ast$)} is a local minimizer of the Kohn-Sham energy functional $E(U)$, we assume in addition that the Hessian of the Kohn-Sham energy functional is bounded from both above and below in a neighborhood of $U^\ast$, that is to say, the following assumption holds.

\begin{assume}\label{Hess_bound}
	There exists $\delta_1>0$, such that for all $ [U] \in B([U^*],\delta_1)$,
	\begin{eqnarray}\label{pos-bnd}
	\nabla^2_GE(U)[D, D] &\ge& \underline{c}\| D \|^2, \forall \ D \in \mathcal{T}_{[U]}\mathcal{G}^N_{N_g},%\notag
	\end{eqnarray}
	and
	\begin{eqnarray}\label{upp-bnd}
	\nabla^2_GE(U)[D] &\le& \bar{c}\| D \|, \forall \ D \in \mathcal{T}_{[U]}\mathcal{G}^N_{N_g},%\notag
	\end{eqnarray}
	%where $\nu_1$ and $\nu_2 > 0$ are constants, and
	where $U^\ast$ is the local minimizer of $E(U)$ and $\bar{c}\ge\underline{c}>0$ are some constants. Here, $B([U],\delta)$ is defined as
	\begin{equation*}
	B([U],\delta) := \{[V]\in \mathcal{G}_{N_g}^N: \min_{P\in\mathcal{O}^{N\times N}}\|U-VP\| \le \delta\}.
	\end{equation*}
\end{assume}
\rv{We mention that the positiveness condition \eqref{pos-bnd} has been justified and used in, e.g., \cite{DLZZ, DWZ, schneider09}, which is related to the spectral gap of Hamiltonian. Meanwhile, the boundedness condition \eqref{upp-bnd} is quite natural as it holds true for any fixed $[U]$ with some constants $\bar{c}_{[U]}>0$, and we just require that there is a uniform upper bound $\bar{c}$ for all $\{\bar{c}_{[U]}\}_{[U] \in B([U^*],\delta_1)}$.}

\begin{remark}
	Under Assumption \ref{Hess_bound}, there exists a positive constant $C$ such that 
	$$\|\nabla_G E(U)\|\le C, ~\forall U\in\mathcal{M}_{N_g}^N,$$
	where $C$ can be chosen as $\sqrt{2}\bar{c}N$.
\end{remark}

Besides, we review some preliminaries on the Grassmann manifold which will be used in the following analysis. 
Let $[U], [W] \in \mathcal{G}_{N_g}^N$, with
$U,W\in \mathcal{M}_{N_g}^N$. We obtain from Lemma A.1 in \cite{DLZZ} that there exists a geodesic
\begin{equation}\label{geodesic1}
\Gamma(t)=[UA\cos{(\Theta t)}A^T+A_2\sin{(\Theta t)}A^T], t\in[0,1],
\end{equation}
such that
\begin{eqnarray*}
	\Gamma(0)=[U], \Gamma(1)=[W].
\end{eqnarray*}
Here, $U^TW=A\cos{\Theta}B^T$ and $W-U(U^TW)=A_2\sin{\Theta}B^T$ is the SVD of $U^TW$ and $W-U(U^TW)$, respectively, $$\Theta=\textup{diag}(\theta_1, \theta_2, \cdots, \theta_N)$$ is a diagonal matrix with $\theta_i\in[0, \pi/2]$ and $$\sin{(\Theta t)}=\textup{diag}(\sin(\theta_1 t), \sin(\theta_2 t), \cdots, \sin(\theta_N t))$$ with a similar notation for $\cos{(\Theta t)}$. Note that $A_2\in \mathcal{M}_{N_g}^N$.

\begin{remark}\label{uni-geo}
	For any $U\in\mathcal{M}^N, D\in\mathcal{T}_{[U]}\mathcal{G}^N$, let $D=ASB^T$ be the SVD of $D$ where $A\in\mathcal{T}_{[U]}\mathcal{G}_{N_g}^N$, $S,B\in\mathbb{R}^{N\times N}$, then there exists an unique geodesic
	\begin{equation}\label{geodesic3}
	\Gamma(t)=[UB\cos{(S t)}B^T+A\sin{(S t)}B^T],
	\iffalse=\textup{exp}_{[U]}(tD)\fi
	\end{equation}
	which start from $[U]$ and with direction $D$ \cite{EAS}. The above expression \eqref{geodesic1} is just a special case with direction $D=A_2\Theta A^T$.
\end{remark}

More specifically, we use macro $[\textup{exp}_{[U]}(t D)]$
to denote the geodesic on $\mathcal{G}_{N_g}^N$ which starts with $[U]$ and with the initial direction $D\in\mathcal{T}_{[U]}\mathcal{G}_{N_g}^N$. We now define the parallel mapping which maps a tangent vector along the geodesic \cite{EAS}.
\begin{definition}
	The parallel mapping $\tau_{\scriptscriptstyle{(U,D,t)}}: \ \mathcal{T}_{[U]}\mathcal{G}_{N_g}^N\to\mathcal{T}_{[\textup{exp}_{[U]}(tD)]}\mathcal{G}_{N_g}^N$ along the geodesic $[\textup{exp}_{[U]}(tD)]$ is defined as
	\begin{equation*}
	\tau_{\scriptscriptstyle{(U,D,t)}}\tilde{D}=\big((-U\sin{(St)}+A\cos{(St)}A^T+(I_N-AA^T)\big)\tilde{D},
	\end{equation*}
	where $D=ASB^T$ is the SVD of $D$.
\end{definition}

It can be verified that \iffalse$\tau_{\scriptscriptstyle{(U,D,t)}}\tilde{D}\in\mathcal{T}_{[\textup{exp}_{[U]}(tD)]}\mathcal{G}_{N_g}^N$ and\fi
\begin{eqnarray}\label{par_invariant}
\|\tau_{\scriptscriptstyle{(U,D,t)}}\tilde{D}\|=\|\tilde{D}\|, \forall \tilde{D}\in\mathcal{T}_{[U]}\mathcal{G}_{N_g}^N.
\end{eqnarray}

To state our theory, we introduce two distances on the Grassmann manifold $\mathcal{G}_{N_g}^N$:  
\begin{equation}\label{dist}
\begin{split}
&\textup{dist}_{cF}([U], [W]) = \min_{P\in \mathcal{O}^{N\times N}} \| U - WP \|, \\
&\textup{dist}_{geo}([U], [W]) =  \| A_2\Theta A^T \|.
\end{split}
\end{equation}

\begin{remark}
	It can be calculated that \cite{EAS}
	\begin{equation*}
	\begin{split}
	&\textup{dist}_{cF}([U], [W]) = \|2\sin{\frac{\Theta}{2}}\|, \\
	&\textup{dist}_{geo}([U], [W]) = \| \Theta \|,
	\end{split}
	\end{equation*}
	which indicate that these two kinds of distance are equivalent. More specifically,
	$$\textup{dist}_{cF}([U], [W])\leq\textup{dist}_{geo}([U], [W])\leq2\textup{dist}_{cF}([U], [W]).$$
	In addition, we see that
	\begin{equation}\label{DeTheta}
	\|D\|  = \|A_2\Theta A^T\| = \|\Theta\|_F = \textup{dist}_{geo}([U], [W]),
	\end{equation}
	where D is the initial direction of the geodesic \eqref{geodesic1}.
\end{remark}

Furthermore, we need the following conclusion, which can be obtained from Remark 3.2 and Remark 4.2 of \cite{Smith94}.
\begin{proposition}\label{Taylor}
	Suppose $E(U)$ is of second order differentiable, then for all $U\in\mathcal{M}_{N_g}^N$, $D\in\mathcal{T}_{[U]}\mathcal{G}_{N_g}^N$, there exists a $\xi\in(0,t)$ such that
	\begin{eqnarray*}
	E(\textup{exp}_{[U]}(tD)) &=& E(U)+t\langle\nabla_GE(\textup{exp}_{[U]}(\xi D)),\tau_{\scriptscriptstyle{(U,D,\xi)}} D\rangle,\label{taylor1} \\
	&=& E(U)+t\langle\nabla_GE(U),D\rangle  \notag\\
	&&+ \frac{t^2}{2}\nabla_G^2E(U)[D, D] + o(t^2\|D\| ^2). \notag 
	\end{eqnarray*}
	and
	\begin{eqnarray*}
	\label{taylor3}
	~~~~\tau_{\scriptscriptstyle{(U,D,t)}}^{-1}\nabla_GE(\textup{exp}_{[U]}(tD)) &=& \nabla_G E(U)+t\tau_{\scriptscriptstyle{(U,D,\xi)}}^{-1}\nabla_G^2E(\textup{exp}_{[U]}(\xi D))[\tau_{\scriptscriptstyle{(U,D,\xi)}} D].
	\end{eqnarray*}
\end{proposition}

Now we are ready to have the local convergence rate of the numerical approximations $\{U_n\}_{n=0}^\infty$ as stated in the following theorem.

\begin{theorem}\label{con-rate}
	Let Assumptions \ref{Assume-energy-decrease} and \ref{Hess_bound} hold true and \rv{assume that} there exists a $\tau>0$ such that $\Delta t_n>\tau, \forall n\in\mathbb{N}_0$. Then for the sequence $\{U_n\}_{n\in\mathbb{N}_0}$ produced by Algorithm \ref{A1} with initial guess $[U_0]\in B([U^*],\delta_1)\subset\mathcal{G}_{N_g}^N$, there exists a constant $\nu\in(0,1)$ such that
	$$E(U_{n+1})-E(U^\ast)\le \nu \big(E(U_{n})-E(U^\ast)\big),$$
	and hence, there exist $C_1, C_2>0$, such that
	$$E(U_n)-E(U^\ast)\le C_1\nu^n \textup{dist}_{geo}([U_0],[U^\ast])^2,$$
	and
	$$\textup{dist}_{geo}([U_n],[U^\ast])\le C_2(\sqrt{\nu})^n \textup{dist}_{geo}([U_0],[U^\ast]).$$
\end{theorem}

\begin{proof}
	For simplicity, we denote $d_n = \textup{dist}_{geo}([U_n],[U^\ast])$. We see that
	\begin{eqnarray*}
		E(U_{n+1})-E(U^\ast) &=& E(U_{n+1}) - E(U_n) + E(U_n) - E(U^\ast) \\
		&\le& -\eta\Delta t_n\|\nabla_G E(U_n)\|^2 +  E(U_n) - E(U^\ast).
	\end{eqnarray*}
	For $U_n$ and $U^\ast \in\mathcal{M}_{N_g}^N$, there exists \rv{a} unique geodesic $\exp_{[U^\ast]}(tD_n)$ such that $$\exp_{[U^\ast]}({\bf0}) = [U^\ast] \ \textup{and} \ \exp_{[U^\ast]}({D_n}) = [U_n],$$
	where ${\bf 0}$ is the zero element on the tangent space $\mathcal{T}_{[U^\ast]}\mathcal{G}_{N_g}^N$. Furthermore, there holds $\|D_n\|=d_n$.
	
	We obtain from \eqref{taylor1} that there \rv{exist} $\xi_{n,1}\in(0,1), \xi_{n,2}\in(0,1)$ such that
	\begin{eqnarray}\label{ener-diff}
	\underline{c}d_n^2\le E(U_n) - E(U^\ast) = \nabla_G^2 E(\exp_{[U^\ast]}({\xi_{n,1}D_n}))[\tau_{\xi_{n,1}}D_n, \tau_{\xi_{n,1}}D_n] \le \bar{c}d_n^2,
	\end{eqnarray}
	and
	\begin{eqnarray}\label{grad-ge-dist}
	\|\nabla_G E(U_n)\| &=& \|\tau_{\xi_{n,2}}^{-1}\nabla_G^2 E(\exp_{[U^\ast]}({\xi_{n,2}D_n}))[\tau_{\xi_{n,2}}D_n]\| \notag\\
	&=&\|\nabla_G^2 E(\exp_{[U^\ast]}({\xi_{n,2}D_n}))[\tau_{\xi_{n,2}}D_n]\|\ge\underline{c}d_n.
	\end{eqnarray}
	Combining \eqref{ener-diff} and \eqref{grad-ge-dist}, we see that
	\begin{equation}\label{grad-bound}
	\|\nabla_GE(U_n)\|^2\ge\frac{\underline{c}^2}{\bar{c}}\big(E(U_n) - E(U^\ast)\big).
	\end{equation}
	Hence, we obtain from the fact that $\{\Delta t_n\}_{n=0}^\infty$ is bounded from below that
	\begin{eqnarray*}
		E(U_{n+1})-E(U^\ast)&\le& -\eta\Delta t_n\|\nabla_G E(U_n)\|^2 +  E(U_n) - E(U^\ast)\\
		&\le& (1-\eta\tau\frac{\underline{c}^2}{\bar{c}})(E(U_n) - E(U^\ast)).
	\end{eqnarray*}
	Finally, we complete the proof by using \eqref{ener-diff} again, and choosing $$\nu = 1-\eta\tau\frac{\underline{c}^2}{\bar{c}}, ~C_1 = \bar{c}, ~C_2 = \sqrt{\frac{\bar{c}}{\underline{c}}}.$$
	 \end{proof}

\section{Some typical orthogonality preserving schemes}\label{sec-ii}
In the previous section, we propose and analyze a  general framework for interpolation based orthogonality preserving schemes for discretizing the gradient flow based Kohn-Sham DFT model \eqref{Gra-Flow}. In that framework, how to determine  the specific form of the auxiliary mapping  $U^{\textup{Aux}}$  and  the time step size $\Delta t_n$  are not given. 
The specific form and the efficiency of Algoirthm \ref{A1}   depend strongly on the definition of the auxiliary mapping  $U^{\textup{Aux}}$ and the choice of the time step size. For example, if the time step size $\Delta t_n$ is chosen to be too large, Assumption \ref{Assume-energy-decrease} may not hold, which may lead to divergence. On the contrary, Theorem \ref{con-rate} indicates that tiny step sizes will lead to slow convergence. In this section,  we will provide some choices for the auxiliary mapping $U^{\textup{Aux}}$, and propose an adaptive approach  for determining the time step sizes. We will prove that our approach can produce time step sizes which can not only satisfy Assumption  \ref{Assume-energy-decrease} but also avoid slow convergence.

\subsection{Auxiliary mapping $U^{\textup{Aux}}$}
We see from the previous discussion that Algorithm \ref{A1} gives a general framework of orthogonality preserving numerical schemes for solving \eqref{generate}, which provides at least a subsequence that converges to the equilibrium point under some mild assumptions. All our analysis \rv{in Section \ref{sec-scheme}} is independent of the specific form of $U^{\textup{Aux}}$ at each interval $[t_n, t_{n+1})$. 
 However, the choice of the auxiliary mapping is one of the keys when we carry out Algorithm \ref{A1}. 
 Here, we provide some potential choices for \rv{the} auxiliary mapping $U^{\textup{Aux}}(t)$.

\newtheorem{exam}{Choice}  
\renewcommand{\labelitemi}{$\bullet$}  

\begin{exam}\label{Linear-Comb}
	$U^{\textup{Aux}}(t) = (1-\alpha_n)U_n + \alpha_n\tilde{U}(t)$, $\alpha_n\in[0, 1]$, $t\in[t_n,t_{n+1})$.
\end{exam}
If we use Crank-Nicolson's strategy \cite{CN47}, which is a widely used second order scheme in time, to discretize \eqref{Gra-Flow}, then we have 
\begin{equation}\label{CN}
U_{n+1} = \big(I_N+\frac{\Delta t}{2}\mathcal{A}_{U_{n+1}}\big)^{-1}\big(I_N-\frac{\Delta t}{2}\mathcal{A}_{U_n}\big)U_n.
\end{equation}
However, it may not preserve the orthogonality of orbitals. Notice that if we choose $\alpha_n = 0$ in Example \ref{Linear-Comb}, then the orbitals can be updated by
\begin{equation}\label{CNL}
U_{n+1} = \big(I_N+\frac{\Delta t}{2}\mathcal{A}_{U_n}\big)^{-1}\big(I_N-\frac{\Delta t}{2}\mathcal{A}_{U_n}\big)U_n,
\end{equation}
which preserves the orthogonality automatically and is an approximation of Crack-Nicolson scheme \eqref{CN} simply by substituting $\mathcal{A}_{U_{n+1}}$ with $\mathcal{A}_{U_n}$ in \eqref{CN}. Hence, we may denote the auxiliary mapping in this case as $U^{\textup{Aux}}_{\textup{CN}}(t) = U_n$. 
Besides, if $\alpha_n$ is chosen to be $1/2$, then we have $U^{\textup{Aux}}(t) = \big(\tilde{U}(t)+U_n\big)/2$. We can see that in this case, the updating formula is the same as the midpoint scheme studied in \cite{DWZ}. Hence, we denote
$$U^{\textup{Aux}}_{\textup{Mid}}(t) = \big(\tilde{U}(t)+U_n\big)/2.$$
Therefore, the framework that we proposed (Algorithm \ref{A1}) contains both the Crank-Nicolson like scheme \eqref{CNL} and the midpoint scheme \cite{DWZ}.

Under the classification mentioned in this paper, we see that the midpoint scheme is implicit and \rv{can not be easily carried out}. Instead, Dai et al. also proposed an explicit approximation to \rv{the} midpoint scheme based on the Picard iteration \cite{DWZ}. More precisely, the midpoint scheme can be replaced approximately by the iterative formulae
$$U^{m}_{n+1/2}(t) = \big(I+\frac{t-t_n}{2}\mathcal{A}_{U^{m-1}_{n+1/2}(\Delta t_n)}\big)^{-1}U_n, ~m=1,2,\dots,$$
where $U^{0}_{n+1/2} = U_n$. This gives us the following choice.

 \begin{exam}\label{appro-mid}
	$U^{\textup{Aux}}(t) = U^m_{n+1/2}(t) \rv{=:} U^{\textup{Aux}}_{\textup{aMid-m}}(t), ~t\in [t_n, t_{n+1}), \forall m\in\mathbb{N}_0$. 
\end{exam}

It is easy to check that $U^{\textup{Aux}}_{\textup{aMid-m}}(t_n) = U_n$. Replacing the midpoint $\frac{U_n+\tilde{U}(t)}{2}$ by \rv{the} approximated midpoints $U^{\textup{Aux}}_{\textup{aMid-m}}(t)$ in the midpoint scheme, we obtain a set of explicit \rv{schemes} and name them as approximated midpoint schemes. They are also included in our interpolation based schemes.

We can of course construct some simpler explicit schemes, e.g., the following one.

 \begin{exam}
	Let 
	\begin{equation}\label{Exp1}
	U^{\textup{Aux}}(t) = U_n - m_n(t-t_n)\nabla_G E(U_n), ~t\in[t_n,t_{n+1})
	\end{equation}
	or 
	\begin{equation}\label{Exp2}
	U^{\textup{Aux}}(t)=2(I+m_n(t-t_n)\mathcal{A}_{U_n})^{-1}U_n-U_n, ~t\in[t_n,t_{n+1}),
	\end{equation}
	where $m_n$ can be arbitrary real number.
 \end{exam}

There are also many other explicit schemes and we will not go further into them. With regard to implicit schemes, we propose the following example motivated by the Verlet algorithm \cite{verlet}. Consider the first order Taylor expansion  of $U(t)$ at $t_{n+\frac{1}{2}} = \frac{t_n+t_{n+1}}{2}$, that is,
\begin{eqnarray*}
	 	U_{n+1} \approx U(t_{n+1/2}+\Delta t_n/2) \approx U(t_{n+1/2}) + \frac{\Delta t_n}{2} \dot{U}(t_{n+1/2}), \label{T1} \\
	U_n \approx U(t_{n+1/2}-\Delta t_n/2) \approx U(t_{n+1/2}) - \frac{\Delta t_n}{2} \dot{U}(t_{n+1/2}). \label{T2}
\end{eqnarray*}
It can be observed that the midpoint scheme uses $\frac{U_n+U_{n+1}}{2}$ to approximate $U_{n+1/2}$ with linear accuracy.
We may further consider the second order Taylor expansion of $U(t)$, which is formulated as
\begin{eqnarray*}
	U_{n+1} \approx U(t_{n+1/2}+\Delta t_n/2) \approx U(t_{n+1/2}) + \frac{\Delta t_n}{2} \dot{U}(t_{n+1/2}) + \frac{\Delta t_n^2}{8} \ddot{U}(t_{n+1/2}), \label{T3} \\
	U_n \approx U(t_{n+1/2}-\Delta t_n/2) \approx U(t_{n+1/2}) - \frac{\Delta t_n}{2} \dot{U}(t_{n+1/2}) + \frac{\Delta t_n^2}{8} \ddot{U}(t_{n+1/2}), \label{T4}
\end{eqnarray*}
where
\begin{eqnarray}\label{2ndder}
\ddot{U}(t) &=& -\frac{d}{dt}\nabla_G E(U(t)) \\
&=&(I-U(t)U(t)^T)\nabla^2E(U(t))[\nabla_GE(U(t))] \notag\\
&&+(\nabla_GE(U(t))U(t)^T+U(t)\nabla_GE(U(t))^T)\nabla E(U(t)) \notag\\
&\rv{=:}&G(U(t)),\notag
\end{eqnarray}
based on which a higher order approximation of the midpoint can be obtained.

	\begin{exam}\label{Verlet}
	\begin{equation*}
U^{\textup{Aux}}(t) = \frac{\tilde{U}(t) + U_n}{2} - \frac{(t-t_n)^2}{8} G(\frac{\tilde{U}(t) + U_n}{2}) =:U^{\textup{Aux}}_{\textup{Verlet}}(t), ~t\in [t_n, t_{n+1}). 
	\end{equation*}
	Here, the operator $G$ can be defined as \eqref{2ndder} or it can be chosen as some approximations of  \eqref{2ndder}.
\end{exam}

 We should emphasize that there are many different choices for the auxiliary mapping $U^{\textup{Aux}}(t)$. Each of them will result in a specific orthogonality preserving scheme for the discretized Kohn-Sham model. The difference lies in   the efficiency, which will be further studied in our future work.

\subsection{Time step sizes}\label{sec-tem}

The choice of the time step sizes is of great importance in the discretization of time dependent problems,  on which many  studies have been done in literature (see, e.g., \cite{qiao-zhang-tang-2011,liao-tang-zhou-2020}). 
 In the numerical analysis for Algorithm 1 provided in Section 3, we require that the time step sizes satisfy Assumption \ref{Assume-energy-decrease}.  
Here,   
we provide a detailed method to help us judge whether or not a preset step size $\Delta t_n$ satifies the energy decrease property of the gradient flow based model \eqref{non-increase}, whose  key  idea is to use the second-order Taylor expansion  to approximate the energy functional $E(\tilde{U}(t_n+\Delta t))$. Note that \rv{a} similar idea has been used in \cite{DZZ}.

By using the second-order Taylor expansion, we have the following  approximation:
  {\small 
  	\begin{equation}\label{eqn-2ndappro}
  	E(\tilde{U}(t_n+\Delta t)) \approx E(\tilde{U}(t_n)) + \Delta t\langle\nabla_G(E(\tilde{U}(t_n))),\tilde{U}'(t_n)\rangle + \frac{\Delta t^2}{2}\nabla_G^2 E(\tilde{U}(t_n))[\tilde{U}'(t_n),\tilde{U}'(t_n)].
  	\end{equation}}
  From the definition of $\tilde{U}$, we see that 
  $$\tilde{U}(t_n)=U_n, ~\tilde{U}'(t_n)=-\nabla_G E(U_n),$$
  and thus \eqref{eqn-2ndappro} becomes
  {\small  
  	\begin{equation}\label{eqn-2ndappro-2}
  	E(\tilde{U}(t_n+\Delta t)) \approx E(U_n) - \Delta t\|\nabla_G(E(U_n)\|^2 + \frac{\Delta t^2}{2}\nabla_G^2 E(U_n)[\nabla_G(E(U_n)),\nabla_G(E(U_n)].
  	\end{equation}}
  Inserting  (\ref{eqn-2ndappro-2}) into (\ref{back-cond}) in Assumption \ref{Assume-energy-decrease}, we have the following inequality for the time step size $\Delta t$
  \begin{eqnarray*}
  	\frac{\|\nabla_G E(U_n)\|^2 - \frac{\Delta t}{2}\nabla_G^2 E(U_n)[\nabla_G E(U_n), \nabla_G E(U_n)]}{\|\nabla_G E(U_n)\|^2} \geq \eta. 
  	\end{eqnarray*}
  Therefore, for a given $\Delta t$ at the  $n$-th iteration, we define the following indicator
  \begin{equation}\label{indicator4}
  \zeta_n(\Delta t) = \frac{\|\nabla_G E(U_n)\|^2 - \frac{\Delta t}{2}\nabla_G^2 E(U_n)[\nabla_G E(U_n), \nabla_G E(U_n)]}{\|\nabla_G E(U_n)\|^2}
  \end{equation}
   to tell us if it is a good step size. 
  If  $\zeta_n(\Delta t)\ge\eta$, we consider $\Delta t$  as a satisfactory time step and  accept it.  
  Otherwise, we instead choose $\Delta t_n$ to be the approximated minimizer of $E(\tilde{U}(t_n+\Delta t))$ with respect to $\Delta t$ at the $n$-th iteration. That is,  
   we choose $$\Delta t_n = \min\{\frac{\|\nabla_GE(U_n)\|^2}{\nabla_G^2E(U_n)[\nabla_GE(U_n), \nabla_GE(U_n)]}, \frac{\theta_n}{\|\nabla_GE(U_n)\|}\},$$ 
which is the minimizer of the right hand side of \eqref{eqn-2ndappro-2} in a small neighbourhood of $0$ to be our final step size.
 
 In summary, we obtain an adaptive strategy to get the time step sizes which will be proved to satisfy Assumption \ref{Assume-energy-decrease}. With this strategy,  the corresponding interpolation based scheme \rv{reads} as the following Algorithm \ref{A-Adap},  where $\delta t_{\min}$ and $\delta t_{\max}$ are the preseted bound \rv{for} the initial step sizes.

\begin{algorithm}\label{A-Adap}
	\caption{Interpolation based scheme with adaptive step sizes}
	%\begin{algorithmic}
	Given $\epsilon, \delta t_{\min}, \delta t_{\max}>0$, $\eta\in(0,1/2)$, initial data $U_0\in\mathcal{M}_{N_g}^N$,
	calculate the gradient $\nabla_G E(U_0)$ and set $n = 0$, $\rv{t_0=0}$\;
	\While{$\| \nabla_G E(U_n) \|> \epsilon$}{
		
		Choose $\theta_n\in(0,1)$ and the initial guess $\Delta t_n^{initial}\in[\delta t_{\min},\delta t_{\max}]$ by some specific strategy;
		
		 Let $\Delta t_n = \Delta t_n^{initial}$;
		
	\If {$\zeta_n(\Delta t_n)<\eta$ or $\Delta t_n > \frac{\theta_n}{\|\nabla_G E(U_n)\|}$
		}{
			
			$\Delta t_n =\min\{\frac{\|\nabla_GE(U_n)\|^2}{\nabla_G^2E(U_n)[\nabla_GE(U_n), \nabla_GE(U_n)]}, \frac{\theta_n}{\|\nabla_GE(U_n)}\|\} $;
		}
		
		Set $t_{n+1} = t_n+\Delta t_n$;
		
		Define $U^{\textup{Aux}}(t)$ on the interval $[t_n, t_{n+1})$ such that $U^{\textup{Aux}}(t_n)=U_n$;
		
		Update $U_{n+1} = \lim_{t\to t_{n+1}^-}\tilde{U}(t)$ with $\tilde{U}(t)$ satisfying \eqref{inner};
		
		Let $n=n+1$, calculate the gradient $\nabla_G E(U_n)$\;
	}
	
\end{algorithm}

For Algorithm \ref{A-Adap},  we have the following theorem, which shows the convergence of our interpolation based scheme with adaptive step sizes.
 
\begin{theorem}\label{Adap-bound}
	If Assumption \ref{Hess_bound} holds and the initial guess $[U_0]\in B([U^*],\delta_1)\subset\mathcal{G}_{N_g}^N$, then there exist $\{\theta_n\}_{n\in\mathbb{N}_0}$ such that for the sequence $\{U_n\}_{n\in\mathbb{N}_0}$ generated by Algorithm \ref{A-Adap}, there holds either $\nabla_G E(U_n)=0$ for some $n\in\mathbb{N}_0$ or 
	\begin{equation}\label{gradtozero}
	\liminf_{n\to\infty}\|\nabla_G E(U_n)\|=0.
	\end{equation}
	Furthermore, there also holds that 
	$$\liminf_{n\to\infty}\textup{dist}_{geo}(U_n,U^\ast) = 0.$$
\end{theorem}

\begin{proof}
	To simplify the notation, we denote $D_n = -\nabla_G E(U_n)$. We see that the time step size $\Delta t_n$ given by Algorithm \ref{A-Adap} should satisfy $\Delta t_n\|D_n\|\le\theta_n$ and
	$$\Delta t_n\|D_n\|^2+\frac{\Delta t_n^2}{2}\nabla_G^2 E(U_n)[D_n,D_n]\le\eta \Delta t_n\|D_n\|^2,~ \forall n\in\mathbb{N}_0.$$
	Define
	\begin{eqnarray*}
		\theta_n&=&\sup\{\tilde{\theta}_n: E(\tilde{U}(t_n+\Delta t))-E(U_n)-\Delta t\|D_n\|^2 \\
		&&-\frac{\Delta t^2}{2}\nabla_G^2 E(U_n)[D_n,D_n]\le-\frac{\eta \Delta t\|D_n\|^2}{2}, \forall \Delta t\le\frac{\tilde{\theta}_n}{\|D_n\|}\}\ge0.
	\end{eqnarray*}
	Then, we obtain from the definition of $E(U_{n+1})$ and $\theta_n$ that
	$$E(U_{n+1})-E(U_n)\le\frac{\eta}{2} \Delta t_n\|D_n\|^2, \forall n\in\mathbb{N}_0,$$
	i.e., \eqref{back-cond} holds.
	
	As for the condition \eqref{sumunbound}, we see that $\Delta t_n$ has only three possible values, that is, $$\Delta t_n=\max{(t_n^{\text{initial}}, \delta t_{\min})},$$ $$\Delta t_n=\frac{\|D_n\|^2}{\nabla_G^2 E(U_n)[D_n, D_n]},$$ or $$\Delta t_n=\frac{\theta_n}{\|D_n\|}.$$ So there is at least one infinite subsequence of $\{n_j\}_{j=0}^{\infty}$, which is, with out loss of generality, also denoted by $\{n\}_{n=0}^{\infty}$, such that
	\vskip 3pt
	
	{\bf Case 1.} $\Delta t_n=\max{(\Delta t_n^{\text{initial}}, \delta t_{\min})}$. We have immediately that $$\sum_{n=0}^{\infty} \Delta t_n\ge\sum_{j=0}^{\infty} \Delta \delta t_{\min} =+\infty.$$
	
	{\bf Case 2.} $\Delta t_n = \frac{\|D_n\|^2}{\nabla_G^2 E(U_n)[D_n, D_n]}$. We obtain from Assumption \ref{Hess_bound} that $\Delta t_n\ge\frac{1}{\bar{c}}$ and
	hence $$\sum_{n=0}^{\infty} \Delta t_n=+\infty.$$
	
	{\bf Case 3.} $\Delta t_n = \frac{\theta_n}{\|D_n\|}$.
	If
	$$\liminf_{n\to\infty}\Delta  t_n>0,$$
	then \eqref{sumunbound} is satisfied. Otherwise, there exists a subsequence of $\{\Delta t_n\}_{n\in\mathbb{N}_0}$, which is also denoted by $\{\Delta t_n\}_{n\in\mathbb{N}_0}$\rv{,} such that
	$\lim_{n\to\infty}\Delta  t_n=0.$
	This simply leads to $\lim_{n\to\infty}\theta_n=0$ since $\|D_n\|$ is bounded from above.
	
	We have that for all $n\in\mathbb{N}_0$, there hold
	\begin{eqnarray*}
		&&E(\tilde{U}(t_n+\Delta t))-E(U_n)-\Delta t\|D_n\|^2-\frac{\Delta t^2}{2}\nabla_G^2 E(U_n)[D_n,D_n]  \\
		&=& E(\tilde{U}(t_n+\Delta t)) - E(\textup{exp}_{[U_n]}(\Delta tD_n)) + E(\textup{exp}_{[U_n]}(\Delta tD_n)) - E(U_n) \\
		&&+\Delta t\|D_n\|^2-\frac{\Delta t^2}{2}\nabla_G^2 E(U_n)[D_n,D_n] \notag =: T_n^{(1)}+T_n^{(2)}, \notag
	\end{eqnarray*}
	where $$T_n^{(1)} = E(\tilde{U}(t_n+\Delta t)) - E(\textup{exp}_{[U_n]}(\Delta tD_n))$$ and $$T_n^{(2)} = E(\textup{exp}_{[U_n]}(\Delta tD_n)) - E(U_n)+\Delta t\|D_n\|^2-\frac{\Delta t^2}{2}\nabla_G^2 E(U_n)[D_n,D_n].$$
	We see from Remark \ref{uni-geo} that there exists a geodesic $[\textup{exp}_{[\tilde{U}(t_n+\Delta t)]}(t\hat{D})]$ such that
	\begin{eqnarray*}
		\textup{exp}_{[\tilde{U}(t_n+\Delta t)]}({\bf 0}) = \tilde{U}(t_n+\Delta t), \ [\textup{exp}_{[\tilde{U}(t_n+\Delta t)]}(\hat{D})] = [\textup{exp}_{[U_n]}(\Delta tD_n)],
	\end{eqnarray*}
	and obtain from \eqref{taylor1} that
	\begin{eqnarray*}
		|T_n^{(1)}| &=& |E(\textup{exp}_{[\tilde{U}(t_n+\Delta t)]}(0\hat{D})) - E(\textup{exp}_{[\tilde{U}(t_n+\Delta t)]}(\hat{D}))| \\
		&=& |\langle\nabla_G E(\textup{exp}_{[\tilde{U}(t_n+\Delta t)]}(\xi \hat{D})), \tau_{\scriptscriptstyle{(\tilde{U}(t_n+\Delta t),\hat{D},\xi)}} \hat{D}\rangle| \\
		&\leq& \|\nabla_G E(\textup{exp}_{[\tilde{U}(t_n+\Delta t)]}(\xi \hat{D}))\|\|\tau_{\scriptscriptstyle{(\tilde{U}(t_n+\Delta t),\hat{D},\xi)}} \hat{D}\| \\
		&\leq& C\|\hat{D}\|,  
	\end{eqnarray*}
	where Assumption \ref{Hess_bound} and \eqref{par_invariant} are used in the last inequality. 
	
	From \eqref{DeTheta}, we have that
	\begin{eqnarray*}
		\|\hat{D}\| &=& \textup{dist}_{geo}([\tilde{U}(t_n+\Delta t)], [\textup{exp}_{[U_n]}(\Delta tD_n)]) \\ &\leq& 2\textup{dist}_{cF}([\tilde{U}(t_n+\Delta t)], [\textup{exp}_{[U_n]}(\Delta tD_n)]) \\
		&\leq& 2\|\tilde{U}(t_n+\Delta t) - \textup{exp}_{[U_n]}(\Delta tD_n)\| \\
		&\leq& 2\big(\|\tilde{U}(t_n+\Delta t) - U_n - \Delta tD_n\| \\
		&&+ \|\textup{exp}_{[U_n]}(\Delta tD_n) - U_n - \Delta tD_n\|\big). 
\end{eqnarray*}
	Notice that $\tilde{U}(t)$ and $\textup{exp}_{[U_n]}(tD_n)$ satisfy
	$$\tilde{U}(t_n)=U_n, ~\tilde{U}'(t_n)=D_n$$
	and
	$$\textup{exp}_{[U_n]}({\bf 0})=U_n, ~{\textup{exp}_{[U_n]}}'({\bf 0}) = D_n,$$
	we have
	 $\|\hat{D}\| = o(\Delta t).$ 
	If the sequence $\{\|D_n\|\}_{n\in\mathbb{N}_0}$ is not bound from below, then we complete the proof. Otherwise, we obtain
	\begin{equation}\label{I1}
	T_n^{(1)} = o(\Delta t\|D_n\|)
	\end{equation}
	since $\|D_n\|$ is bounded from both above and below.
	As for $T_n^{(2)}$, \eqref{taylor1} gives that
	\begin{equation}\label{I2}
	T_n^{(2)} = o(\Delta t^2\|D_n\|^2).
	\end{equation}
	Combining \eqref{I1} and \eqref{I2}, we arrive at
	\begin{eqnarray*}
		&&E(\tilde{U}(t_n+\Delta t))-E(U_n)-\Delta t\|D_n\|^2-\frac{\Delta t^2}{2}\nabla_G^2 E(U_n)[D_n,D_n]  \\
		&=& T_n^{(1)}+T_n^{(2)} = o(\Delta t\|D_n\|), ~~ \forall n\in\mathbb{N}_0.
	\end{eqnarray*}
	
	Note that the definition of $\theta_n$ implies that for all $n$, there exists a $$\Delta t_n^\ast\in(\frac{\theta_n}{\|D_n\|}, \frac{\theta_n+\frac{1}{n}}{\|D_n\|}),$$ such that
	\begin{eqnarray}\label{contra}
	o(\Delta t_n^\ast\|D_n\|)  &=& E(\tilde{U}(t_n+\Delta t_n^\ast))-E(U_n) \notag\\
	&&+\Delta t_n^\ast\|D_n\|^2-\frac{{\Delta t_n^\ast}^2}{2}\nabla_G^2 E(U_n)[D_n,D_n] \notag\\
	&>&\frac{\eta \Delta t_n^\ast\|D_n\|^2}{2}.% \notag\\
	\end{eqnarray}
	Hence, it is easy to see that
	$$0\le\lim_{n\to\infty}t_n^\ast\|D_n\|\le\lim_{n\to\infty}\big(\theta_n+\frac{1}{n}\big)=0.$$
	Finally, we obtain by letting $n\to\infty$ in \eqref{contra} that
	$$0\ge\lim_{n\to\infty}\frac{\eta}{2}\|D_n\|,$$
	which together with \eqref{grad-ge-dist} completes the proof.
\end{proof}

\rv{
\section{Numerical experiments}\label{sec-ne}
\renewcommand*{\thefootnote}{\arabic{footnote}}
In this section, we apply one of our proposed schemes \rv{to solve} the discretized Kohn-Sham DFT model for some typical systems, \rv{including} benzene$(C_6H_6)$, aspirin($C_9H_8O_5$) and Fullerin($C_{60}$), to validate our theoretical results. More specifically, we test the scheme (Algorithm \ref{A-Adap}) with auxiliary mapping \eqref{Exp1} and with $m_n$ being chosen as $1/2$. All of our experiments are carried out on LSSC-IV cluster and the coding is built based on the software package Octopus\footnote[1]{Octopus: {\color{blue}{octopus-code.org/wiki/Main\_Page}}} (Version 4.0.1). Among all our experiments, we set $\eta=1$e$-4$, $\epsilon = 1$e$-12$, $\delta t_{\min}=1$e$-20$, and $\Delta t_n^{initial} = 0.1$, $\theta_n = 0.8$, for all $n$. Here and hereafter, we denote this specific scheme as GF-EX scheme.

We first test the orthogonality preserving property of our scheme. \rv{To this end, we define the orthogonality violation of the iterative orbital $U_n$ at the $n$-th iteration as
\[
\varepsilon_n = \| U_n^T U_n- I_N \|_F
\]
and show the curves for $\{\varepsilon_n\}_n$ in Figure \ref{f1}, of which the $x$-axis stands for the number of iteration $n$ and the $y$-axis is the value of $\varepsilon_n$.}

\begin{figure} 
\centering
\includegraphics[width=.48\textwidth]{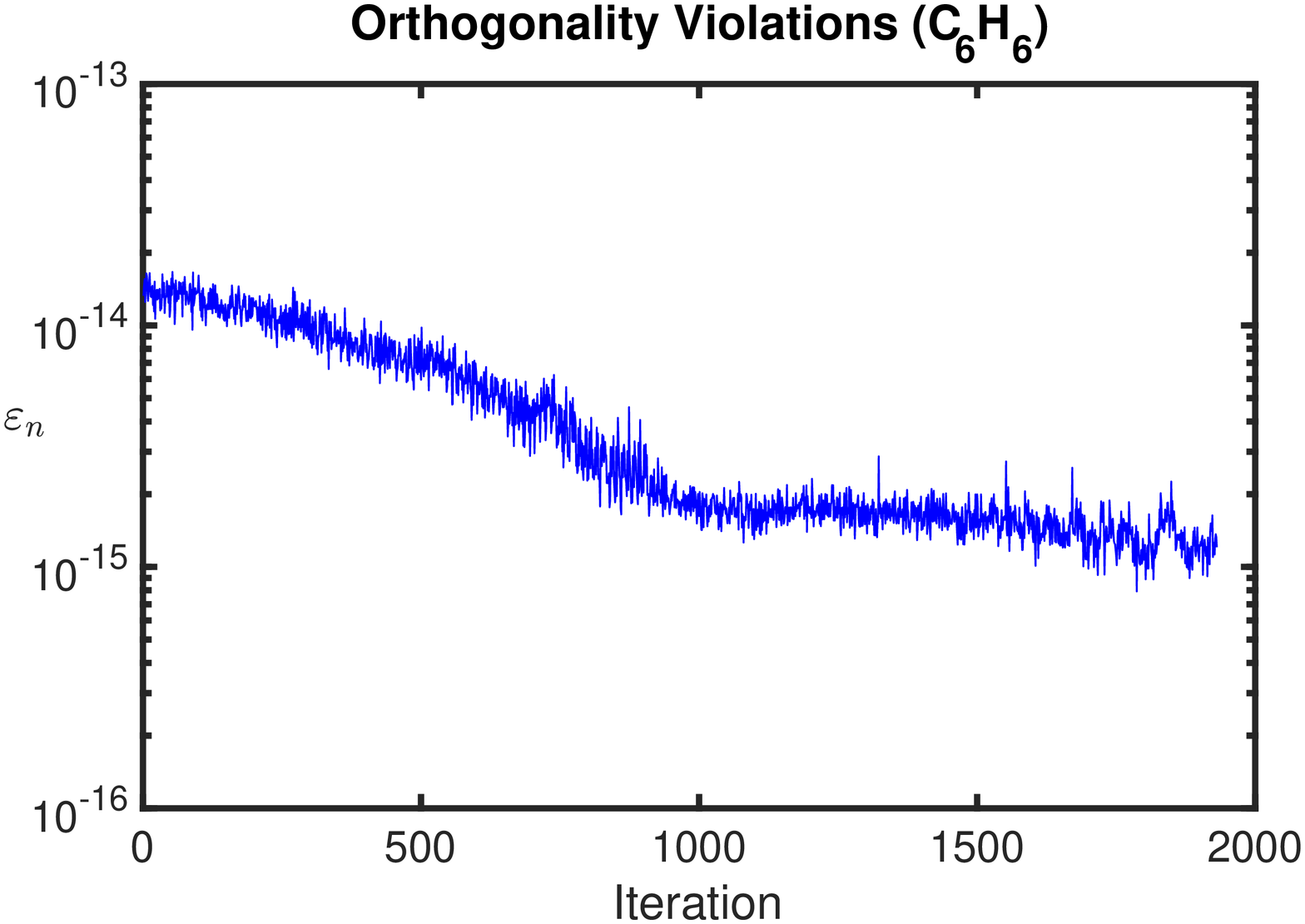} 
\vskip 15pt
\includegraphics[width=.48\textwidth]{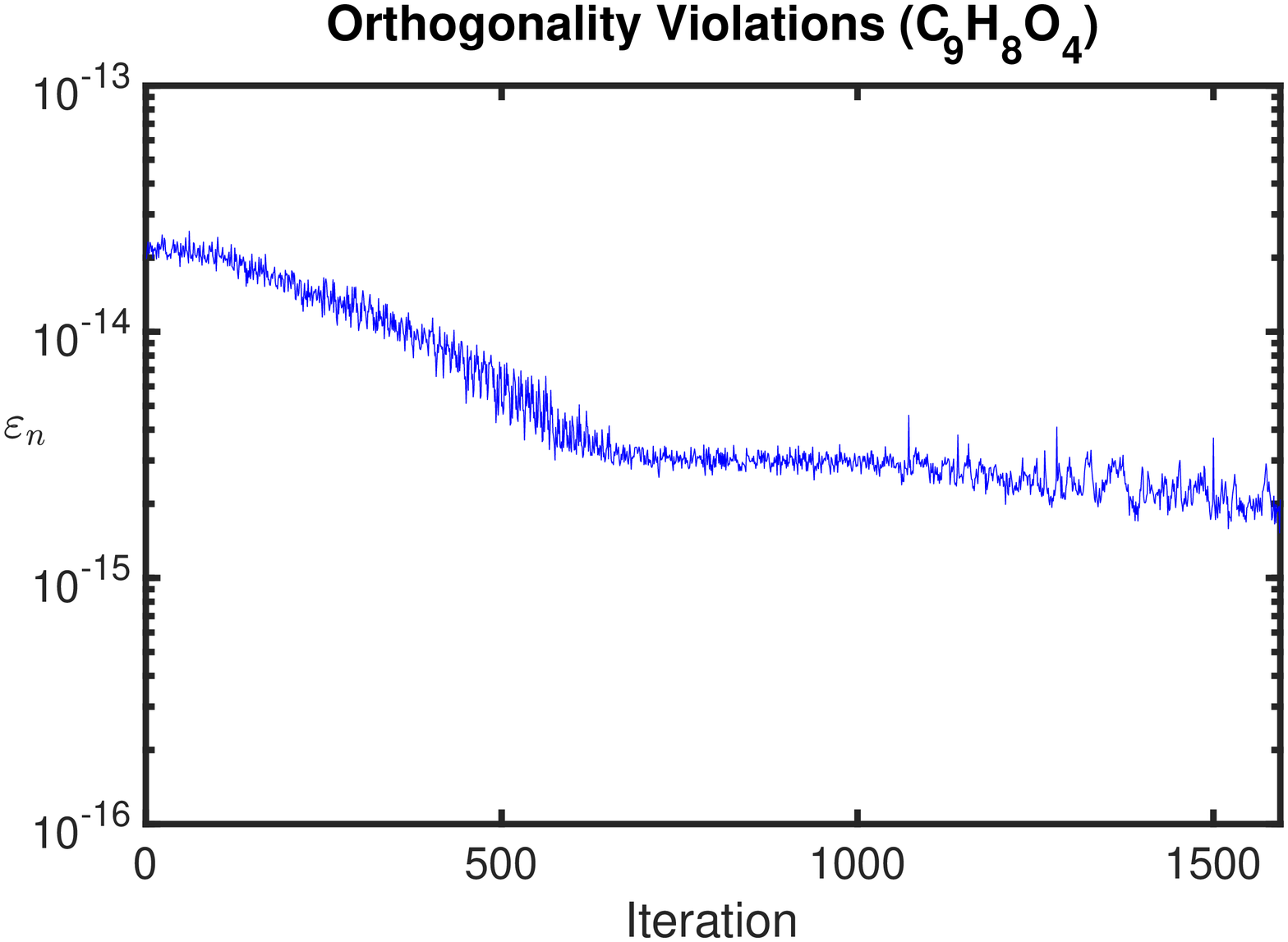} ~~~~
\includegraphics[width=.48\textwidth]{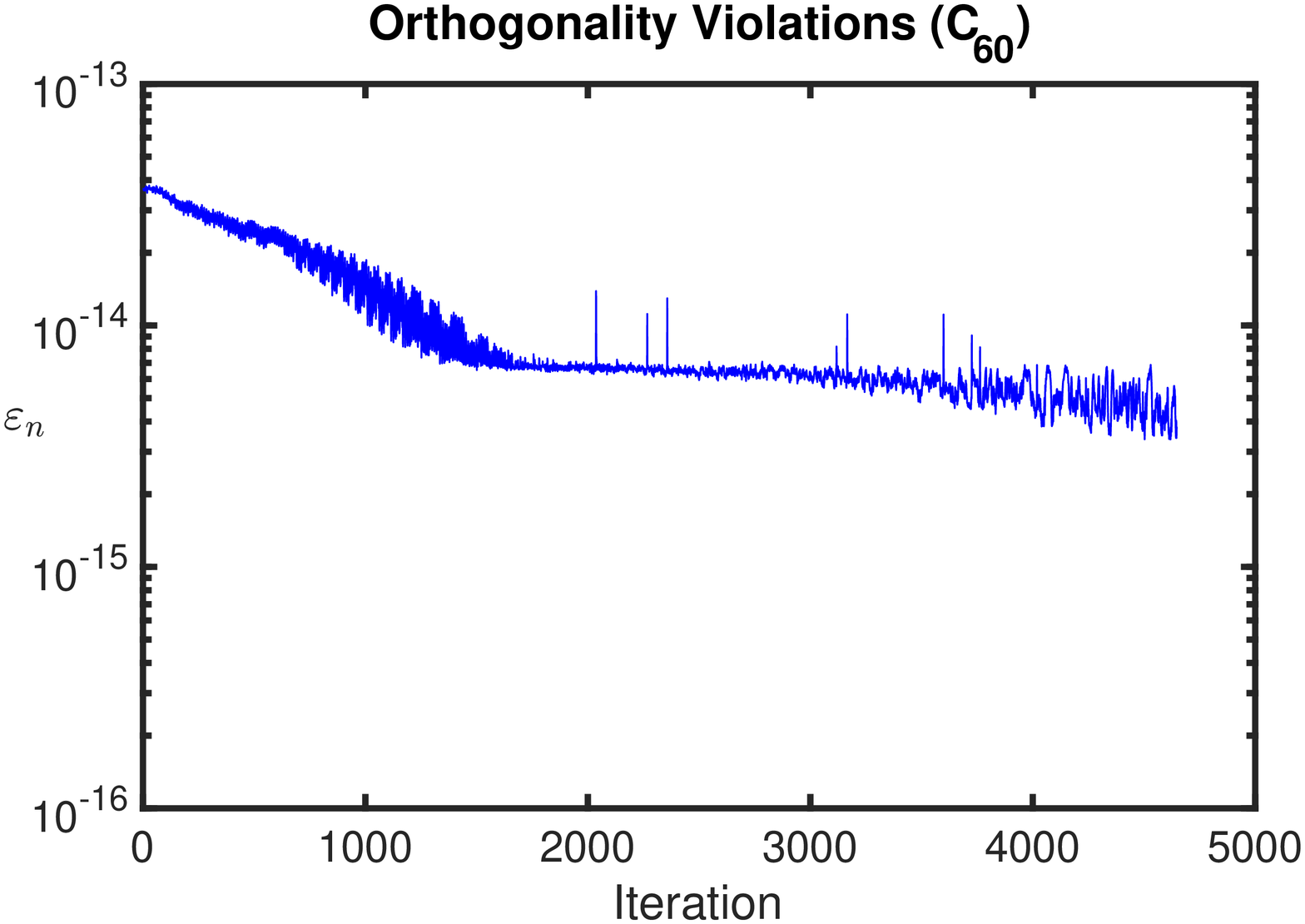} 
\vskip 15pt
\caption{Orthogonality violations obtained by GF-EX for different systems.} 
\label{f1} 
\end{figure} 

It can be observed from Figure \ref{f1} that the orthogonality violations for all tested systems always lie in the interval (1e-15,1e-13) during iterations, which indicates that the GF-EX scheme indeed preserves the orthogonality of iterative orbitals well. 

Then we showcase the convergent states obtained by our scheme in Table \ref{t1}, in which the reference ground state energy $E_{\text{min}}$ is obtained by the SCF iterations on Octopus and $U_{\text{final}}$ stands for the last iterative orbitals when the iteration meets the stopping criterion.

\begin{center} 
\begin{table} [!htbp]
\caption{Numerical results obtained by the scheme GF-EX.} 
\label{t1} 
\begin{center} 
{\small
\begin{tabular} {|c | c| c| c|} 
\hline
System & Reference energy $E_{\text{min}}$ (a.u.) & $E(U_{\text{final}})$ (a.u.) & $\| \nabla_G E(U_{\text{final}})\| $  \\
\hline
Benzene$(C_6H_6)$  & -3.74246025E+01 & -3.74246025E+01  & 9.92E-13  \\
\hline
Aspirin$(C_9H_8O_4)$  & -1.20214764E+02 & -1.20214764E+02  & 6.69E-13 \\
\hline
Fullerin$(C_{60} )$  & -3.42875137E+02 & -3.42875137E+02  & 9.91E-13  \\
\hline
\end{tabular} } 
\end{center} 
\end{table} 
\end{center} 

We see from Table \ref{t1} that our scheme can indeed produce approximations that converge to the ground state.
The following Figures \ref{f2}-\ref{f3} illustrate the convergence curves for the error of energy and the norm of $\nabla_G E$ obtained by our scheme, respectively, which give an intuitive look for the numerical behaviour of the GF-EX scheme. 
\begin{figure} 
\centering
\includegraphics[width=0.48\textwidth]{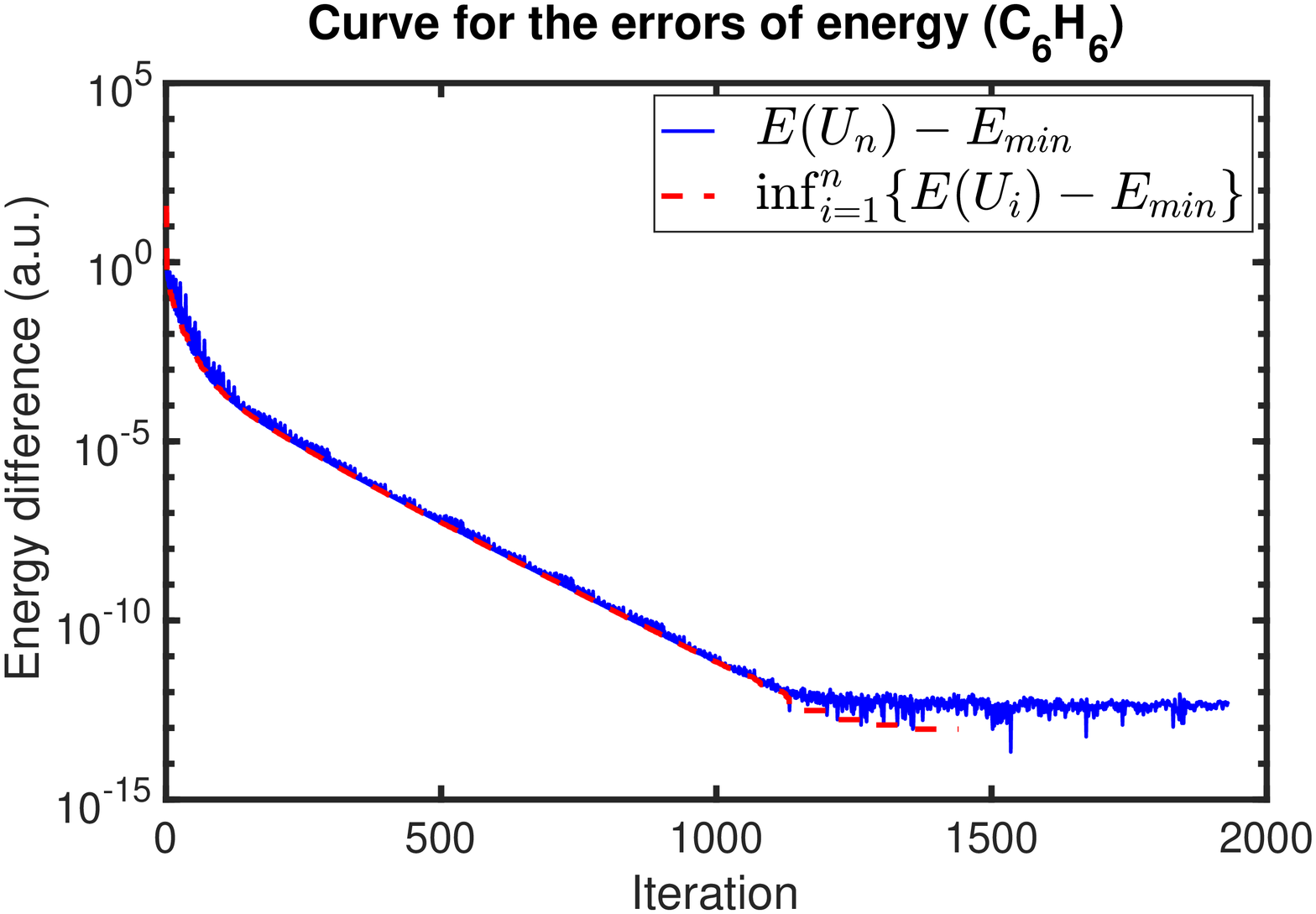} 
\vskip 15pt
\includegraphics[width=0.48\textwidth]{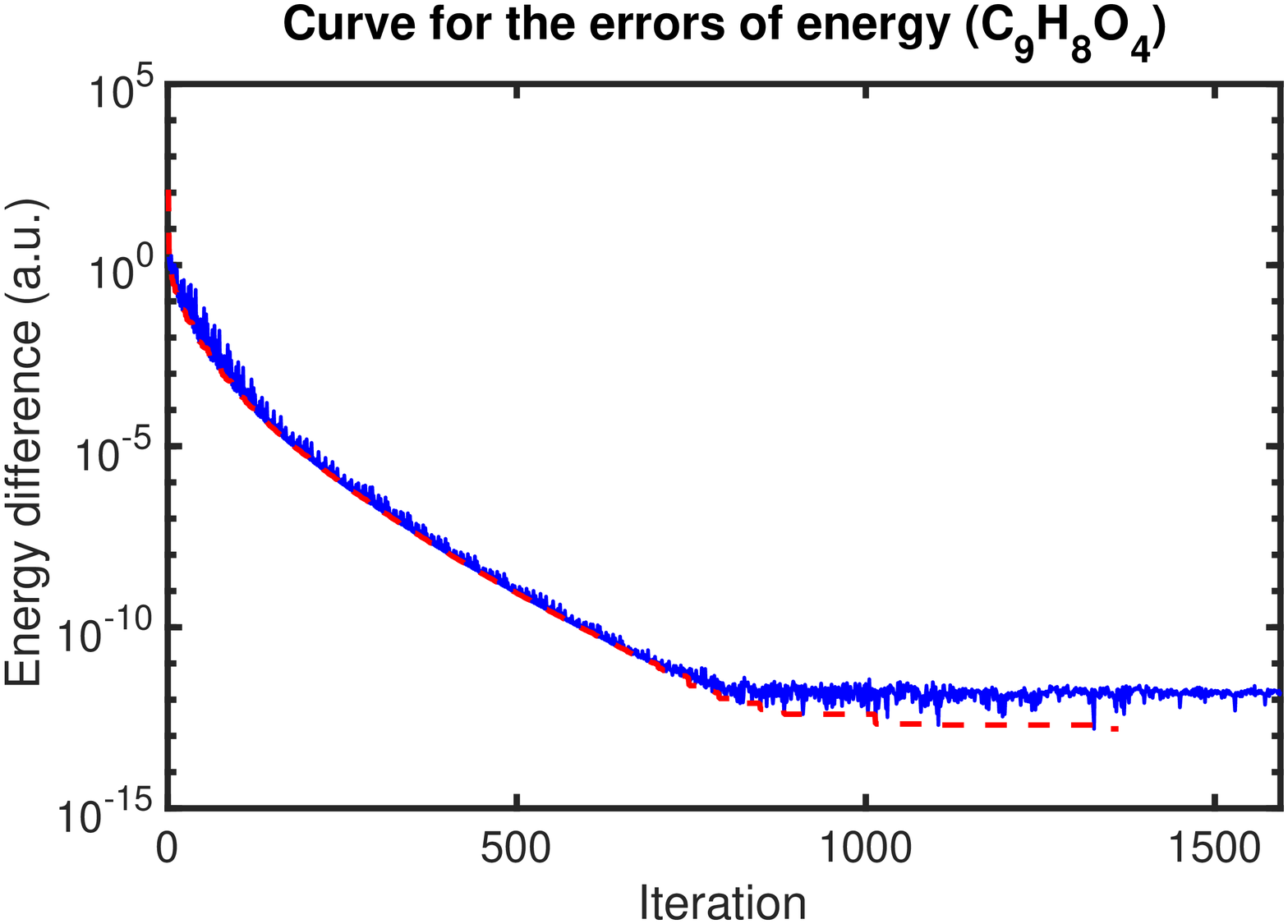} ~~~~
%\vskip 5pt
\includegraphics[width=0.48\textwidth]{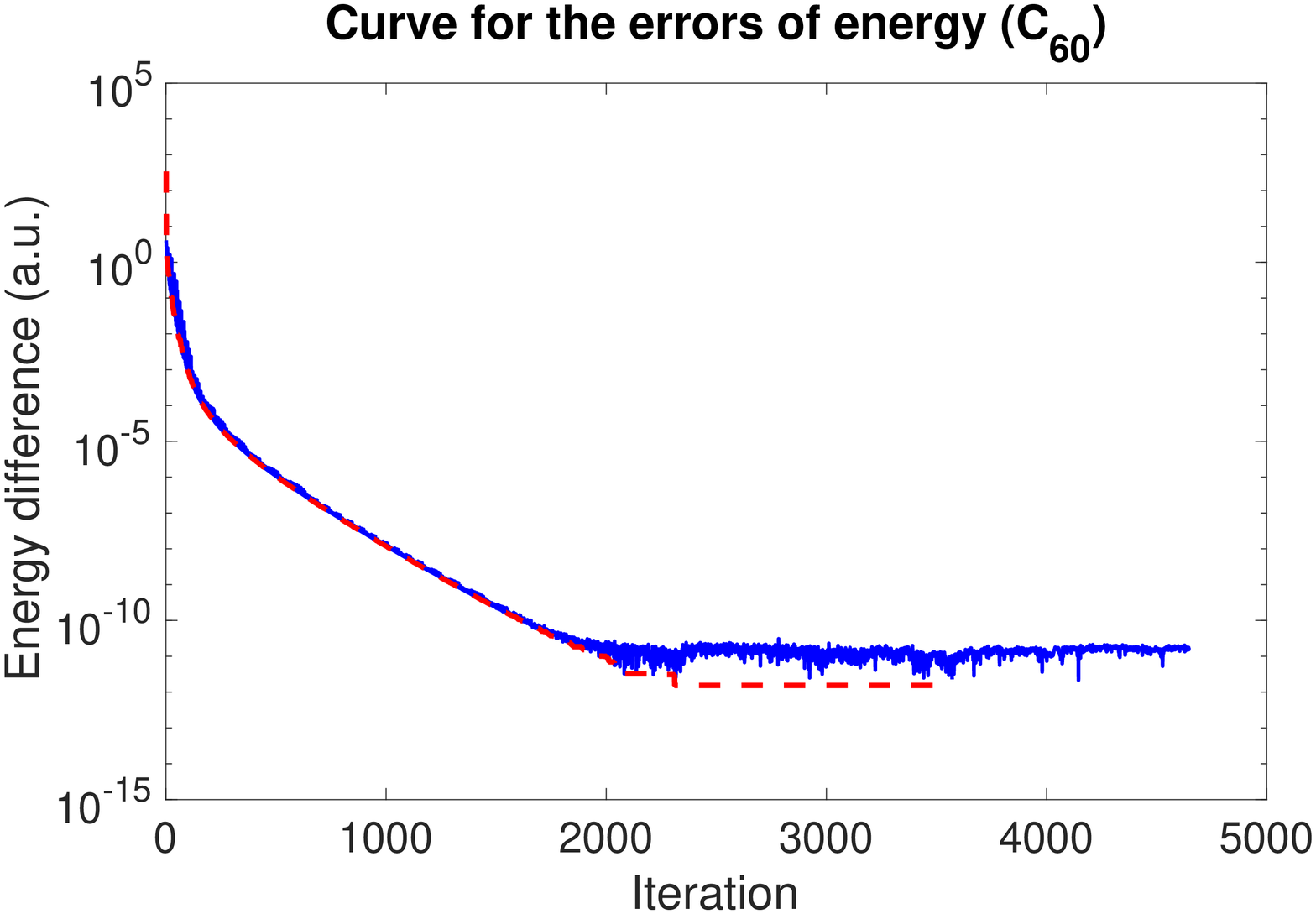} 
\vskip 15pt
\caption{Curves for the error of energy obtained by GF-EX for different systems.} 
\label{f2} 
\end{figure} 

\begin{figure} 
\centering
\includegraphics[width=0.48\textwidth]{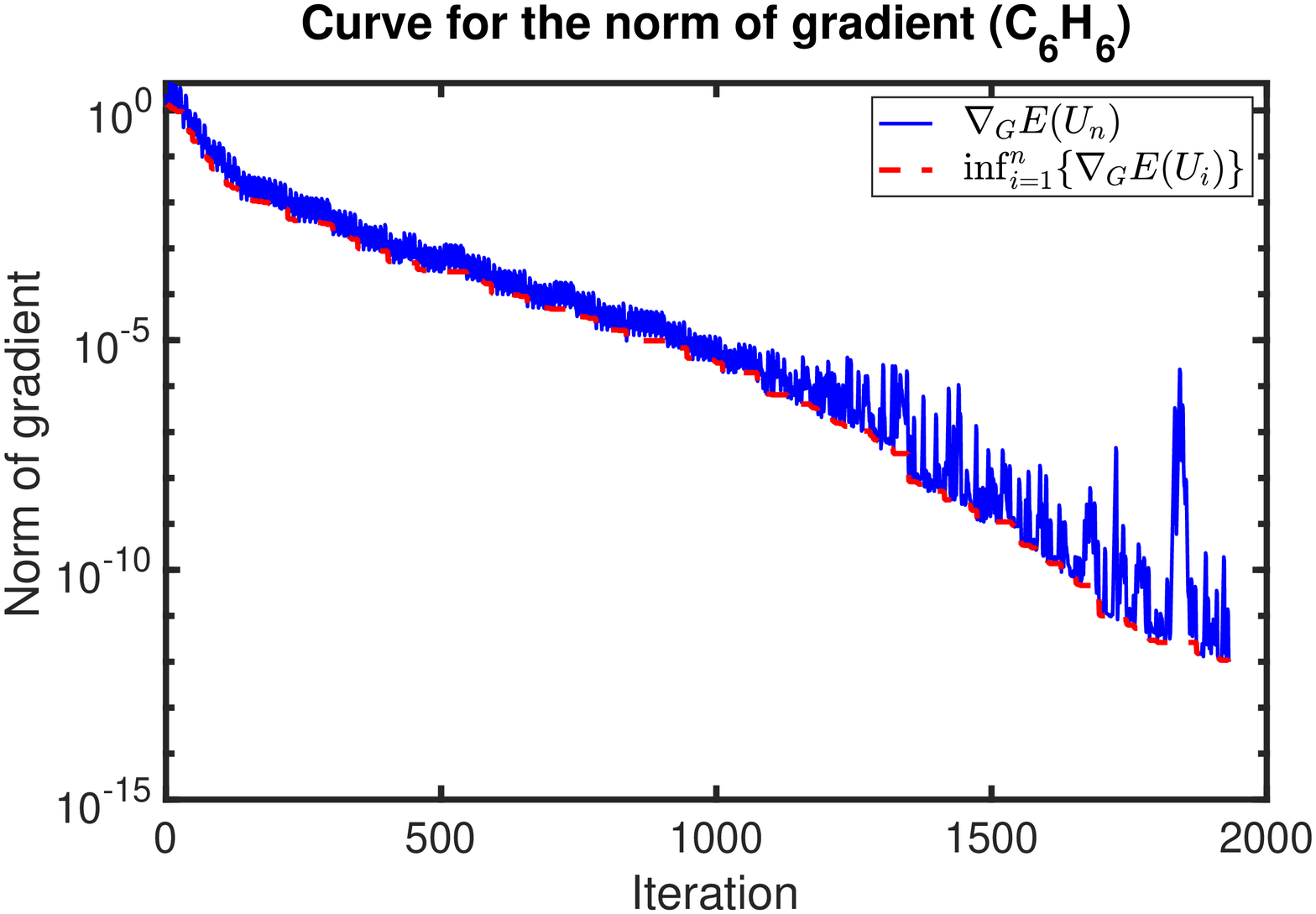} 
\vskip 15pt
\includegraphics[width=0.48\textwidth]{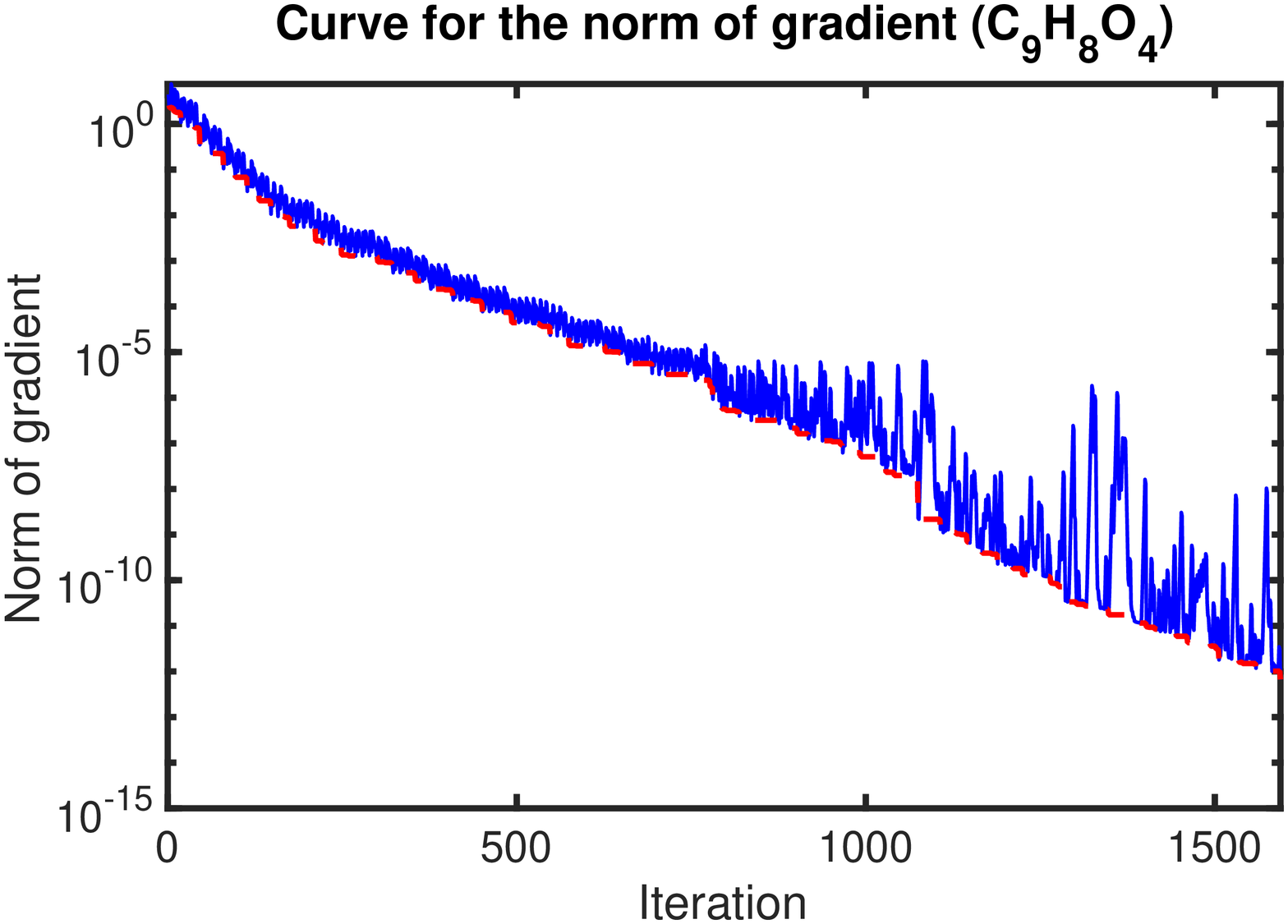} ~~~~
%\vskip 5pt
\includegraphics[width=0.48\textwidth]{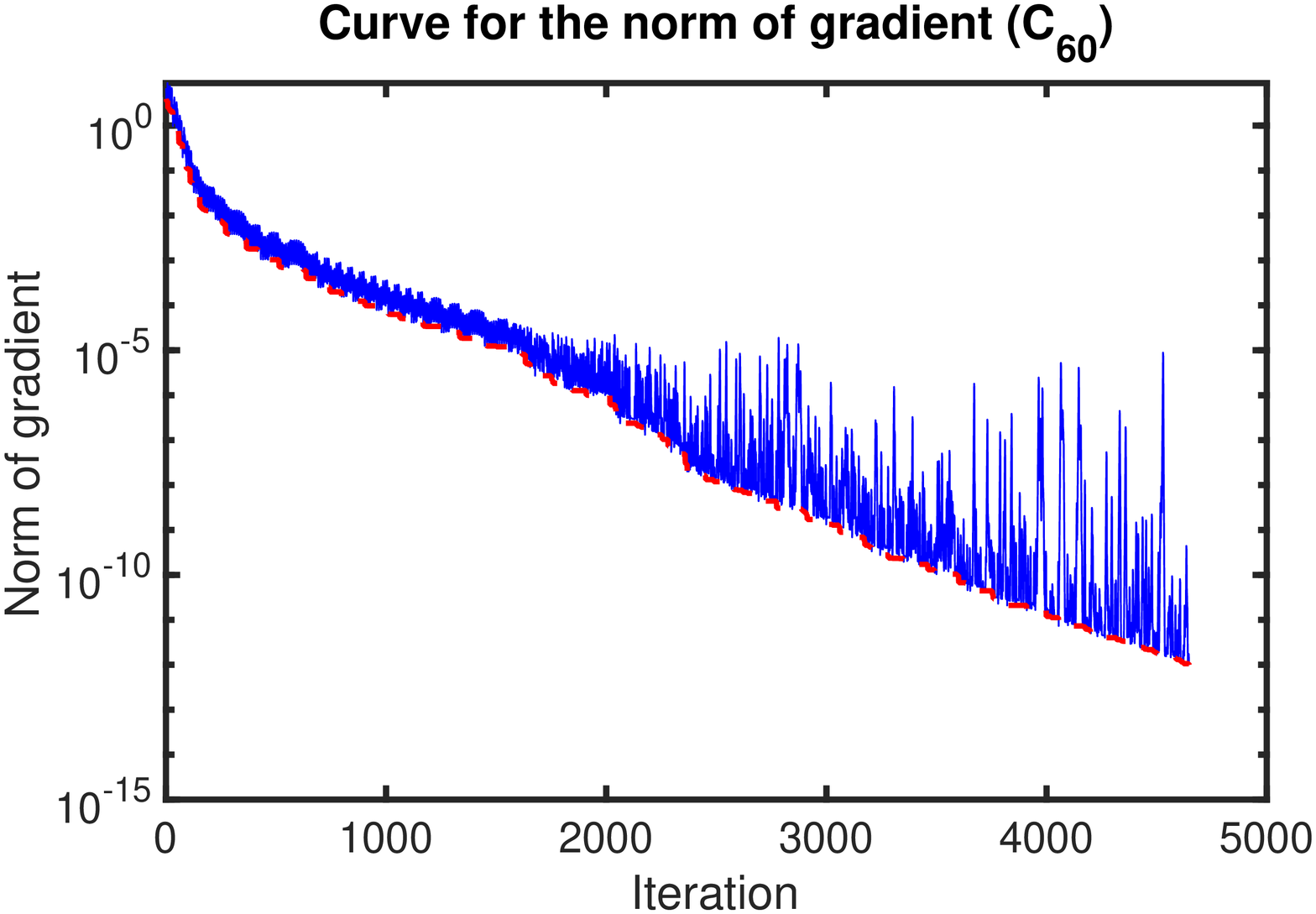} 
\vskip 15pt
\caption{Curves for $\|\nabla_G E\|$ obtained by GF-EX for different systems.} 
\label{f3} 
\end{figure} 

In the energy plots (Figure \ref{f2}), the $y$-axis indicates the energy difference between $E(U_n)$ and $E_{\text{min}}$. Among both Figures \ref{f2} and \ref{f3}, the red lines show the asymptotical infimum of iterations, which are defined as
\[
\inf_n\{E(U_n) - E_{\text{min}}\} = \min_{i\in\{1,2,\cdots,n\}} \{E(U_i) - E_{\text{min}}\},
\]
and
\[
\inf_n \|\nabla_G E(U_n)\| = \min_{i\in\{1,2,\cdots,n\}}\|\nabla_G E(U_i)\|.
\]
We observe that the red line in Figure \ref{f2} terminates earlier than the blue one. The reason is that we have achieved a lower energy than the reference energy $E_{\text{min}}$ during iteration, which makes $\inf_n\{E(U_n)-E_{\min}\}$ negative and thus can no longer be shown in the log scale plots. 

The above two figures show the convergence of both the energy and the norm of gradient clearly, which is consistent to our \rv{theory} 
except that the total energy is not monotonically decreasing. This is because the parameter $\theta_n$ is simply chosen as a fix contant $0.8$, which is usually larger than the one that we defined in the proof of Theorem \ref{Adap-bound}, 
%not chosen overcautiously in our experiments but is chosen as a fix constant $0.8$ 
as we understand that smaller $\theta_n$ may lead to smaller step sizes and thus slow down the convergence. 
}

 \section{Concluding remarks}\label{sec-cln}
In this paper, we proposed and analyzed a general framework of orthogonality preserving schemes for approximating the Kohn-Sham orbitals, from which we can  obtain a class of orthogonality preserving schemes.   We proved the convergence and derived the local \rv{exponential} convergence rate of the framework under some mild and reasonable assumptions. In addition, we provided some typical choices for the auxiliary mapping  which lead to   several  orthogonality preserving schemes.  We then presented an efficient approach to obtain the desired time step sizes that satisfy the assumptions required in our analysis. \rv{Finally, we apply one of the explicit schemes that we proposed as an example to verify our theory. Due to the great flexibility on choosing both auxiliary mapping and step sizes in our framework, we will systematically study, apply and compare the schemes generated by our framework based on numerical experiments on electronic structure calculations \iffalse of some typical systems\fi in our future work.}

\section*{Acknowledgments}
This work was supported by the National Key R $\&$ D Program of China under grants
2019YFA0709600 and 2019YFA0709601, and the National Natural Science Foundation of China under
grants 12021001. The authors would like to thank the anonymous referees for their helpful comments 
and suggestions that enriched the content and improved the presentation of this paper.


\begin{thebibliography}{1}
	
	\bibitem{AbMaSe} {\sc P.-A. Absil, R. Mahony, and R. Sepulchre}, {\em Optimization algorithms
		on matrix manifolds}, Princeton University Press, Princeton, 2008.
	
	\bibitem{anderson65} {\sc D. G. Anderson}, {\em  Iterative procedures for nonlinear integral equations}, J. Assoc. Comput. Mach.,
	12(1965), pp. 547-560.
	
	
	\bibitem{BLL} {\sc Z. Bai, R.-C. Li, and D. Lu}, {\em Optimal convergence rate of self-consistent field iteration for
		solving eigenvector-dependent nonlinear eigenvalue problems}, arXiv:2009.09022(2020).
	
	
	\bibitem{Cances01}{\sc E. Cances}, {\em Self-consistent field algorithms for Kohn-Sham models with fractional occupation
		numbers}, J. Chem. Phys., 114 (2001), pp. 10616-10622.
	
	
	\bibitem{cances12} {\sc E. Cances, R. Chakir, and Y. Maday}, {\em Numerical analysis of the planewave discretization of
		some orbital-free and Kohn-Sham models}, M2AN, 46 (2012), pp. 341-388.
	
	
	\bibitem{Cances21}{\sc E. Cances, G. Kemlin, and A. Levitt}, {\em Convergence analysis of direct minimization and self-consistent 
		iterations}, SIAM J. Matrix Anal. Appl., 42 (2021), pp. 243-274. 
	
	\bibitem{chen-dai14} {\sc 	
		H. Chen, X. Dai, X. Gong, L. He, and A. Zhou}, {\em  Adaptive finite element approximations for Kohn-
		Sham models}, Multiscale Model. Simul., 12(2014), pp. 1828-1869.
	\bibitem{chen13} {\sc 
		H. Chen, X. Gong, L. He, Z. Yang, and A. Zhou}, {\em  Numerical analysis of finite dimensional approximations of Kohn-Sham equations}, Adv. Comput. Math., 38(2013), pp.  225-256.
	
	
	\bibitem{CN47} {\sc J. Crank, and P. Nicolson}, {\em A practical method for numerical evaluation of solutions of partial 
		differential equations of the heat conduction type},
	Proc. Camb. Phil. Soc., 43 (1947), pp. 50-67.
	
	\bibitem{DLZZ} {\sc X. Dai, Z. Liu, L. Zhang, and A. Zhou},
	{\em A conjugate gradient method for electronic structure calculations},
	SIAM J. Sci. Comput., 39 (2017), pp. 2702-2740.
	
	\bibitem{DLZxZ} {\sc X. Dai, Z. Liu, X. Zhang, and A. Zhou},
	{\em A parallel orbital-updating based optimization
		method for electronic structure calculations},
	J. Comput. Phys., 445 (2021), 110622.
	
	\bibitem{DWZ} {\sc X. Dai, Q. Wang, and A. Zhou},
	{\em Gradient flow based discretized Kohn-Sham density
		functional theory },
	Multiscale Model. Simul., 18 (2020), pp. 1621-1663.
	
	\bibitem{DZZ} {\sc X. Dai, L. Zhang, and A. Zhou},
	{\em An adaptive step size strategy for orthogonality constrained line search methods},
	arXiv:1906.02883(2019).
	
	
	\bibitem{dai-zhou15} {\sc X. Dai  and A. Zhou},
	{\em Finite element methods for electronic structure calculations} (in Chinese), SCIENTIA SINICA Chimica, 45 (2015), pp. 800-811.   
 
	\bibitem{DingZhou} {\sc J. Ding and A. Zhou},  
	{\em A spectrum theorem for perturbed bounded linear operators}, 
	Appl. Math. Comput., 201 (2008), pp. 723-728.
	
	\bibitem{EAS} {\sc A. Edelman, T.A. Arias, S.T. Smith}, {\em The geometry of algorithms with
		orthogonality constraints}. SIAM J. Matrix Anal. Appl., 20 (1998), pp. 303-353.
	
	\bibitem{Frey09}{\sc C. Freysoldt, S. Boeck, and J. Neugebauer}, {\em Direct minimization technique for metals in
		density functional theory}, Phys. Rev. B, 79 (2009), 241103.
	
	\bibitem{GLCY} {\sc B. Gao, X. Liu, X. Chen, and Y. Yuan}, {\em A new first-order
		framework for orthogonal constrained optimization problems},
	SIAM J. Optim., 28 (2017), pp. 302-332.
	
	\bibitem{GLY} {\sc B. Gao, X. Liu, and Y. Yuan}, {\em Parallelizable algorithms for optimization problems with orthogonality constraints}, SIAM J. Sci. Comput., 41 (2019), pp. A1949-A1983
	
	\bibitem{HK} {\sc P. Hohenberg and W. Kohn}, {\em Inhomogeneous electron gas},
	Phys. Rev. B., 136 (1964), pp. 864-871.
	
	 
	\bibitem{johnson88} {\sc D.D. Johnson}, {\em  Modified Broyden’s method for accelerating convergence in self-consistent calculations}, Phys. Rev. B, 38(1988),  pp. 12807-12813
	
	
	\bibitem{KS65} {\sc W. Kohn and L. J. Sham}, {\em Self-consistent equations including
		exchange and correlation effects}, Phys. Rev. A., 140 (1965), pp. 4743-4754.
	
	 	\bibitem{liao-tang-zhou-2020}{\sc H. Liao, T. Tang, and T. Zhou}, {\em A second-order and nonuniform
		time-stepping maximum-principle preserving scheme for time-fractional
		Allen-Cahn equations}, J. Comput. Phys., 414 (2020),  109473.
	

	\bibitem{LWWUY} {\sc X. Liu, Z. Wen, X. Wang, M. Ulbrich, and Y. Yuan},
	{\em On the analysis of the discretized Kohn-Sham density functional theory},
	SIAM J. Numer. Anal., 53 (2015), pp. 1758-1785.
	
	\bibitem{LWWY} {\sc X. Liu, X. Wang, Z. Wen, and Y. Yuan}, {\em On the convergence
		of the self-consistent field iteration in Kohn-Sham density functional theory},
	SIAM J. Matrix Anal. Appl., 35 (2014), pp. 546-558.
	
	
	\bibitem{lin-lu-ying2019} {\sc L. Lin L, J. Lu, and L. Ying}, {\em  Numerical methods for Kohn-Sham density functional theory},  Acta Numerica, 28(2019), pp. 405-539.
	
	
	\bibitem{lin-yang13} {\sc L. Lin  and C. Yang}, {\em  Elliptic preconditioner for accelerating the self-consistent field iteration in 
		Kohn-Sham density functional theory}, SIAM J. Sci. Comput., 35(2013), pp. S277-S298.
	
	
	\bibitem{Matin04} {\sc R. Martin}, {\em Electronic Structure: Basic Theory and Practical
		Methods}, Cambridge University Press, London, 2004.
		
	\rv{\bibitem{M3GR}{\sc M. A. L. MARQUES, N. T. MAITRA, F. M. S. MOGUEIRA, E. K. U. GROSS, AND A. RUBIO,} eds., {\em Fundamentals of Time-Dependent Density Functional Theory}, Lecture Notes in Physics, Vol. 837, Heidelberg: Springer Berlin Heidelberg, Berlin, 2012.}
	
	\bibitem{Marzari97}{\sc N. Marzari, D. Vanderbilt, and M. C. Payne}, {\em Ensemble density-functional theory for ab initio molecular dynamics of metals and finite-temperature insulators}, Phys. Rev. Lett.,
	79 (1997), 1337.
	
	
	\bibitem{parr-yang94} {\sc R. G. Parr and W. Yang}, {\em  Density-Functional Theory of Atoms and Molecules}, Oxford University Press, New York, 1994.
	
	
	\bibitem{payne92} {\sc M.C. Payne,   M. P. Teter, D. C. Allen,  T. A. Arias,  and J. D. Joannopoulo},  {\em Iterative minimization techniques for ab initio total energy calculation: Molecular dynamics and conjugate gradients}, Rev. Mod. Phys. 64(1992),  1045-1097.
	
	\bibitem{pulay80} {\sc P. Pulay}, {\em Convergence acceleration of iterative sequences: the case of scf iteration}, Chem. Phys.
	Lett., 73(1980), pp. 393-398.
	
	\bibitem {pulay82} {\sc P. Pulay}, {\em  Improved SCF convergence acceleration}, J. Comput. Chem., 3(1982), pp. 556-560.
	
 
	\bibitem{qiao-zhang-tang-2011} {\sc Z. Qiao, Z. Zhang, and T. Tang}, {\em An adaptive time-stepping strategy
	for the molecular beam epitaxy models}, SIAM J. Sci. Comput., 33(2011), pp. 1395-1414.
	
	\rv{\bibitem{RG1984} {\sc E. Runge and E. K. U. Gross,} {\em Density functional theory for time-dependent systems}, Phys. Rev. Lett., 52 (1984), 997-1000.}

	  
	\bibitem{Saad10}{\sc Y. Saad, J. R. Chelikowshy, and S. M. Shontz}, {\em Numerical methods for electronic structure
		calculations of materials}, SIAM Rev., 52 (2010), pp. 3-54.
	
	\bibitem{schneider09} {\sc R. Schneider, T. Rohwedder, A. Neelov, and J. Blauert},
	{\em Direct minimization for calculating invariant subspaces in density fuctional
		computations of the electronic structure}, J. Comput. Math., 27 (2009), pp. 360-387.
	
	\bibitem{Smith94} {\sc S. T. Smith}, {\em Optimization techniques on Riemannian
		manifolds}, in Fields Institute Communications, Vol. 3, AMS, Providence, RI, 1994, pp. 113-146.
	
	\bibitem{verlet} {\sc L. Verlet}, {\em Computer ‘experiments’ on classical fluids. 1. Thermodynamical properties of LennardJones molecules},
	Physical Review., 159 (1967), pp. 98-103.
	
	 
	\bibitem {yang-wang07} {\sc C. Yang, J.C. Meza, and L. Wang}, {\em A trust region direct constrained minimization algorithm for the
		Kohn-Sham equation},  SIAM J. Sci. Comput., 29(2007), pp. 1854-1875.
	
	
	\bibitem{YGM} {\sc C. Yang, W. Gao, and J. Meza}, {\em On the convergence of the
		self-consistent field iteration for a class of nonlinear eigenvalue problems},
	SIAM J. Matrix Anal. Appl., 30 (2009), pp. 1773-1788.
	
	\rv{\bibitem{YSHH}{\sc L. Yang, Y. Shen, Z. Hu, and G. Hu,} {\em An implicit solver for the time-dependent Kohn-Sham equation.} Numer. Math. Theor. Meth. Appl., 14(2020), 261-284.}
	
	\bibitem{zhang-shen08}	{\sc D. Zhang, L. Shen, A. Zhou, and X. Gong}, {\em Finite element method for solving Kohn-Sham
		equations based on self-adaptive terahedral mesh}, Phy. Lett. A, 372 (2008), pp. 5071-5076.
	
	\bibitem{ZZWZ} {\sc X. Zhang, J. Zhu, Z. Wen, and A. Zhou},
	{\em Gradient type optimization methods for electronic structure calculations},
	SIAM J. Sci. Comput., 36 (2014), pp. 265-289.
	
	\bibitem{ZBJ} {\sc Z. Zhao, Z. Bai, and X. Jin},
	{\em A Riemannian Newton algorithm for nonlinear eigenvalue problems},
	SIAM J. Matrix Anal. Appl., 36 (2015), pp. 752-774. 
	
	\bibitem{zhou-wang-2018}  {\sc Y. Zhou, H. Wang,  Y. Liu,  X. Gao, and H. Song}, {\em 
		Applicability of Kerker preconditioning scheme to the self-consistent density functional theory calculations of inhomogeneous systems}, 
	Phys. Rev. E, 97(2018),  033305.
	
	
\end{thebibliography}
\end{document}